\documentclass{amsart}
\usepackage{amssymb}
\newtheorem{thm}{Theorem}[section]
\newtheorem{cor}[thm]{Corollary}
\newtheorem{lem}[thm]{Lemma}
\newtheorem{prop}[thm]{Proposition}
\theoremstyle{definition}
\newtheorem{defn}[thm]{Definition}
\theoremstyle{remark}
\newtheorem{rem}[thm]{Remark}
\numberwithin{equation}{section}
\newcommand{\To}{\longrightarrow}

\newcommand{\inv}{^{-1}}
\newcommand{\C}{\mathbb C}
\newcommand{\Z}{\mathbb Z}
\newcommand{\R}{\mathbb R}

\newcommand{\K}{\mathbb K}

\newcommand{\Hecke}{\mathcal H} 

\newcommand{\open}{j}
\newcommand{\closed}{i}

\newcommand{\x}{\times}

\newcommand{\Sh}{\operatorname{Sh}}

\newcommand{\si}{\sigma}

\newcommand{\Supp}{\operatorname{Supp}}
\newcommand{\Cosupp}{\operatorname{Cosupp}}
\newcommand{\Hom}{\operatorname{Hom}}
\newcommand{\Rhom}{\operatorname{RHom}}
\newcommand{\RHom}{\operatorname{RHom}}

\newcommand{\cone}{{\operatorname{cone}}}
\newcommand{\Const}{\underline\C} 
\newcommand{\p}{{\underline p}}
\newcommand{\q}{{\underline q}}
\newcommand{\0}{{\underline 0}}
\newcommand{\dsheaf}{{\mathcal D}}
\newcommand{\dcomplex}{{\hat{\mathcal D}^\bullet}}
\newcommand{\ICsheaf}{\mathcal {IC}}
\newcommand{\IpCsheaf}{\mathcal I_\p\mathcal C}

\newcommand{\Vect}{\text{(vector spaces)}}
\newcommand{\pt}{pt}
\newcommand{\pl}{pl}

\begin{document}

\title[An introduction to perverse sheaves]
{An introduction to perverse sheaves}
\author{Konstanze Rietsch}%
\address{King's College, University of London}%
\email{rietsch@mth.kcl.ac.uk}%


\thanks{
The author was funded by the Violette and Samuel Glasstone
Foundation at Oxford throughout most of the preparation of these
lectures, and is currently a Royal Society Dorothy Hodgkin
Research Fellow at King's College, University of London.}


\date{\today}

\maketitle
\section*{Introduction}

The aim of these notes is to give an introduction to perverse
sheaves with applications to representation theory or quantum
groups in mind. The perverse sheaves that come up in these
applications are in some sense extensions of actual sheaves on
particular algebraic varieties arising in Lie theory (for example
the nilpotent cone in a Lie algebra, or a Schubert variety, or
the moduli space of representations of a Dynkin quiver). We will
ultimately therefore be interested in perverse sheaves on
algebraic varieties.

The origin of the theory of perverse sheaves is M.~Goresky and
R.~MacPherson's theory of intersection homology
\cite{GM:IH1,GM:IH2}. This is a purely topological theory, the
original aim of which was to find a topological invariant similar
to cohomology that would carry over some of the nice properties of
homology or cohomology of smooth manifolds also to singular
spaces (especially Poincar\'e duality). While the usual cohomology
of a topological space can be defined sheaf theoretically as
cohomology of the constant sheaf, the intersection homology turns
out to be the cohomology of a certain complex of sheaves,
constructed very elegantly by P.~Deligne. This complex is the main
example of a perverse sheaf.

The greatest part of these notes will be taken up by explaining
Goresky and MacPherson's intersection homology and Deligne's
complex. And throughout all of that our setting will be purely
topological. Also we will always stay over $\C$. The
deepest result in the theory of perverse sheaves on
algebraic varieties, the Decomposition Theorem of Beilinson,
Bernstein, Deligne and Gabber \cite{BBD:PerverseSheaves}, 
will only be stated. And the notes will end with 
an application, the intersection
cohomology interpretation of the Kazhdan-Lusztig polynomials.

The overriding goal of this exposition (written by a non-expert)
is to give a hopefully broadly accessible first introduction to
perverse sheaves. It is intended more to give the flavor and some
orientation without delving too much into technical details. The
hope is that readers wishing to see something in more detail or
greater generality should be able to orient themselves very
quickly in the existing literature to find it. Moreover these
lectures are very far from comprehensive, and the bibliography
mostly just reflects the sources that I happened across and found
useful. Many more pointers in all the different directions of the
theory can be found in the references themselves. For a
particularly nice overview of the many applications of perverse
sheaves the two ICM addresses \cite{Lus:ICM} and \cite{MacPh:ICM}
are highly recommended.

These notes are based on seven lectures given first at Oxford
University in the Spring of 2002, and then repeated in condensed
form at the Fields Institute, as part of the ICRA X conference. I
am grateful to the listeners in both places for their questions
and comments, and to Professor Ringel for instigating these lectures.
Special thanks go to Kevin McGerty and Catharina Stroppel for each
sending me their helpful comments on the written notes.

\section{Lecture - local systems, an introduction}

Let $X$ be a topological space. We will always assume $X$ to be a
nice, sensible topological space (locally compact, Hausdorff,
paracompact, with a countable basis, locally simply connected).
For example picture your favorite simplicial or CW complex,
manifold or real/complex algebraic variety.  We begin with a
quick review of sheaves. Good references are for example
\cite{Iversen:SheafCohBook, KashScha:SheavesBook}.
Any results from
algebraic topology we use can be found in most textbooks,
such as \cite{Dold:AlgTopBook,Munkres:AlgTopBook,Rotman:AlgTopBook}.

\subsection{Sheaves}
Let $\K$ be a field.
A sheaf of $\K$-vector spaces, $\mathcal F$, on $X$ is a
contravariant functor
\begin{equation*}
\begin{array}{ccc}
\left (\begin{array}{c}\text{open sets in $X$}\\
\text{and inclusions}\end{array}\right )&\longrightarrow &\left
(\begin{array}{c}\text{$\K$-vector spaces}\\
\text{and linear maps}\end{array}\right)\\
U&\longmapsto &\mathcal F(U)\\
V\hookrightarrow U&\longmapsto &r_{UV}:\mathcal F(U)\to \mathcal F(V):
s\mapsto s|_V
\end{array}
\end{equation*}
obeying the {\it sheaf axiom~:} For any collection of open sets
$\{V_i\}_{i\in I}$ in $X$ and $s_i\in\mathcal F(V_i)$ that are
compatible with one another in the sense that $s_i|_{V_i\cap
V_j}=s_j|_{V_i\cap V_j}$ for all $i,j\in I$ there exists a unique
$s\in \mathcal F(\bigcup_{i\in I}V_i)$ such that $s|_{V_i}=s_i$
for all $i\in I$. The elements of  $\mathcal F(U)$ are called {\it
sections} of $\mathcal F$ over $U$.

Loosely, the sheaf axiom says that sections are determined by
their restrictions to any open cover. In a way a sheaf
is a device for dealing with properties that are local in
nature and measuring the transition from local to global. We list
some common examples. These are widely spread, but mostly we will
be interested in sheaves related to 4. below.
\begin{enumerate}
\item The continuous $\R$-valued
functions on a topological space $X$
form a sheaf $C^0_X$. So $C^0_X(U)$ is the vector space of continuous
maps  $U\to\R$ and for $U\subset V$,  $r_{UV}$ is the usual restriction map.
\item
Similarly if $X$ is a smooth manifold $C_X^\infty$ is the sheaf of
smooth functions on $X$. A system of differential equations on
$X$ determines a subsheaf of $C_X^\infty$ with sections given by
the local solutions. If $X$ is an analytic manifold (with analytic
gluing maps between the charts) then one also has the sheaf of
analytic functions on $X$. For an algebraic variety there is the
sheaf of regular (algebraic) functions.
\item
A vector bundle $E\to X$ gives rise to a sheaf of local sections. 
For example we have the sheaf of vector fields coming
from the tangent bundle, or the sheaf of $1$-forms associated to
the cotangent bundle.
\item
A central role in topology is played by the {\it constant sheaf~}
$\Const_X$ on $X$, which is defined by
 $$
\Const_X(U):=\{\text{locally constant maps $U\to
\C$}\}.
 $$
So if $U$ has finitely many connected components we have
 $$
\Const_X(U)\cong \C^{\#\{\text{connected components of $U$}\}}.
$$
Completely analogously there is a constant sheaf $\underline V_X$
for any vector space $V$. In particular we have the zero-sheaf,
$\underline 0_X$, with $\underline 0_X(U)=0$ for all $U$.
\end{enumerate}

A morphism $\mathcal F\to \mathcal G$ between sheaves on $X$ is
defined to be a natural transformation of the functors, that is, a
collection of compatible linear maps $\mathcal F(U)\to \mathcal
G(U)$ for all open sets $U$ in $X$.

We denote the category of sheaves of $\C$-vector spaces on $X$ by
$\Sh(X)$. This is an abelian category. The zero object is given by
the zero-sheaf $\underline 0_X$. The direct sum $\mathcal F\oplus
\mathcal G$ of two sheaves is just given in the obvious way by
$(\mathcal F\oplus \mathcal G)(U)=\mathcal F(U)\oplus\mathcal
G(U)$ with the induced restriction maps. Also the kernel of a
morphism $f:\mathcal F\to\mathcal G$ in the category may be
defined by $\mathcal {K}er(f)(U)=\ker(f(U):\mathcal
F(U)\to\mathcal G(U))$. Note that whether a section in $\mathcal
F(U)$ is sent to zero under $f(U)$ is a property that can be
checked locally, so $\mathcal {K}er(f)$ is indeed a sheaf.

Cokernels are a bit more tricky to define. The definition on
individual spaces of sections is not local enough and only gives
a presheaf (a contravariant functor of the same kind but that
does not necessarily obey the sheaf axiom). However there is a
standard way to `sheafify', i.e. turn any presheaf into a sheaf
by forcing the sheaf axiom, to get the correct definition. We'll
skip the details, but see any textbook. In any case, an exact
sequence of sheaves may be characterized by being exact for small
enough neighborhoods, or by being exact on {\it stalks~:}

\begin{defn}Let $x\in X$ and $\mathcal F$ a sheaf (or presheaf) on
$X$.
The stalk of
$\mathcal F$ at $x$ is the $\C$-vector space
 $$
\mathcal F_x=\lim_{\overset{\longrightarrow} {\{U\subset X\
\text{open}\ |\  x\in U\}}} \mathcal F(U).
 $$
In other words an element in $\mathcal F_x$ can be represented by
a pair $(U,s)$ of an open neighborhood $U$ of $X$ and a section $s
\in \mathcal F(U)$. And two such pairs represent the same element
of the stalk if they agree restricted to some small enough neighborhood
of $x$. This element is the {\it germ} of the section $s$
and denoted by $s_x$.

Any map of sheaves $f:\mathcal F\to\mathcal G$ on $X$
induces a map on stalks $f_x:\mathcal F_x\to \mathcal G_x$. For a
section $s\in\mathcal F(U)$ the {\it support} of $s$ is
defined by $\Supp(s)=\{x\in U\ |\ s_x\ne 0 \text{ in $\mathcal F_x$}\}$.
Note, using the definition of $s_x$, that
$\Supp(s)$ is automatically closed(!) inside $U$. 
\end{defn}

A sequence,
 $
 \underline 0_X\to\mathcal E\to\mathcal F\to\mathcal
 G\to\underline 0_X,
 $
in $\Sh(X)$ is exact precisely if
$0\to \mathcal E_x\to\mathcal F_x\to\mathcal G_x\to 0$ is
an exact sequence of vector spaces for all $x\in X$. In particular the
functor from $\Sh(X)$ to vector spaces of
taking stalks is exact.  Also for any
presheaf $\tilde{\mathcal F}$, the sheafification $\mathcal F$ is a sheaf with
$\mathcal F_x=\tilde{\mathcal F_x}$ for all $x\in X$.

\begin{rem}[Examples]\label{r:stalks}
For a sheaf of analytic functions the germs at $x$ can be thought
of as Taylor series around $x$. The stalk of the constant sheaf
$\underline V_X$ at any point $x$ is just the vector space $V$.

Let $x_0\in X$. Then  another example of a sheaf is the
skyscraper sheaf $\mathcal S$ at $x_0$. Define $\mathcal S(U)=\C$
precisely if $x_0\in U$ and otherwise zero, with the obvious
restriction maps. Then $\mathcal S_{x_0}=\C$ while all other
stalks vanish.
\end{rem}

\subsection{A few functors and exactness}\label{s:restrictions}
The stalk functor is a special case of a {\it restriction
functor}. Suppose $\open=\open_Z:Z\hookrightarrow X$ is the
inclusion of a locally closed subset. Then there is an exact
functor $\open^*:\Sh(X)\to\Sh(Z)$ defined by $\open^*(\mathcal
F)(V)={\underset{\To}{\lim}}\,\mathcal F(U)$, where the limit is
taken over $U\supset V$ which are open in $X$. When $Z$ is open
in $X$ the limit is not required. One also writes
$\open^*(\mathcal F)=\mathcal F|_{Z}$. This functor has the
property $\open^*(\mathcal F)_z=\mathcal F_z$ for all $z\in Z$.

Other important functors that we have straight away are
 \begin{align*}
 \Gamma~: \Sh(X)&\to (\text{$\C$-vector spaces})& \Gamma_c~:\Sh(X)
&\to (\text{$\C$-vector spaces})\\
 \mathcal F&\mapsto \Gamma(X,\mathcal F):=\mathcal F(X), &
 \mathcal F&\mapsto \Gamma_c(X,\mathcal
   F),
 \end{align*}
where $\Gamma_c(X,\mathcal F)$ is the space of compactly
supported global sections. For open $U\subset X$ we also have the
functors $\Gamma(U,\ \ )=\Gamma\circ {\open_U}^*$ and
$\Gamma_c(U,\ \ )=\Gamma_c\circ {\open_U}^*$ of sections or
compactly supported sections over $U$. These sections functors
are in general only left exact.

\begin{rem}
\label{r:deRham} The failure of right exactness of $\Gamma$ comes
from the problem that local properties need not be true globally.
For example suppose $X$ is a smooth manifold and $\Omega^i_X$ the
sheaf of smooth $i$-forms on $X$, say with complex coefficients.
Then the complex of sheaves
 $$
\cdots \to 0\to\Const_X\to\Omega^0_X\to\Omega^1_X\to\cdots
 $$
with the usual differential is exact, since this is locally true
by Poincar\'e's Lemma (every closed form on $\R^n$ is exact).
However applying $\Gamma$ to the degree $\ge 0$ half of this
complex we obtain the de Rham complex for computing $H^*(X,\C)$.

In other words think of the cohomology of the de Rham complex as
giving an obstruction for contractibility of a manifold $X$. Then
this example comes down to the simple observation that while any
manifold is locally contractible, it is not necessarily contractible
globally.
\end{rem}

\subsection{Local systems}
\begin{defn} A sheaf $\mathcal L$ on $X$ is called {\it locally constant}
if every $x\in X$ has a neighborhood $U$ such that for all $y\in
U$, the canonical map
\begin{equation*}
\mathcal L(U)\to \mathcal L_y
\end{equation*}
is an isomorphism. A {\it local system} is a locally constant
sheaf with finite dimensional stalks. If $X$ is connected then all
these stalks automatically have the same dimension. This dimension
is called the rank of the local system.
\end{defn}
Here are some examples.
\begin{enumerate}
\item The standard example of a locally constant sheaf is of course
the constant sheaf $\underline V_X$ for a vector space $V$.
\item Suppose $X$ is a connected $n$-manifold. The top degree
cohomology with compact
supports $H^n_c(X,\C)$ detects whether $X$ has a coherent orientation
or not: if so it is one-dimensional, spanned by a so-called
`orientation class', otherwise zero.
The orientation sheaf $\mathbb O_X$ may be defined by
$\mathbb O_X(U)= H^n_c(U,\C)^*$. Note that an open inclusion
$U\hookrightarrow U'$ induces
a map $H^n_c(U,\C)\to H^n_c(U',\C)$ and therefore the dual is required
to define a sheaf.
$\mathbb O_X$ is a rank one local system
which is constant precisely if $X$ is orientable.
\item \label{i:Lalpha} The local solutions to a homogeneous system of differential
equations on a manifold can form a locally constant sheaf. For
example fix $\alpha\in \C$ and take the differential equation,
\begin{equation}\label{alpha}
\frac{df}{dz}-\frac{ \alpha}z f=0,\tag{$*_\alpha$}
\end{equation}
on the punctured plane $X=\C^*$.
Then we may associate to this equation a sheaf $\mathcal L^\alpha$ on
$X$ by
\begin{equation}\label{e:Lalpha}
\mathcal L^\alpha (U)=
\{\text{ complex analytic functions on $U$ satisfying \eqref{alpha}}\}.
\end{equation}
If $U$ is simply connected then any branch of $z^\alpha$ gives a
solution, and these are related to one another by a scalar
multiple. So $\mathcal L^\alpha (U)=\left <z^\alpha\right >_\C$.
We have that $\mathcal L^\alpha$ is a rank one local system.
$\mathcal L^\alpha$ is trivial (the constant sheaf) precisely if
$\alpha\in\Z$. If the differential operator $\frac{d}{dz}-\frac{
\alpha}z $ is replaced by $(\frac{d}{dz}-\frac{ \alpha}z )^m$ one
obtains a rank $m$ local system. This is the starting point of
Riemann-Hilbert correspondence. See for example
Chapter~7 of \cite{Kirwan:IHBook} and then \cite{Borel:DModuleBook}.
\end{enumerate}

\subsection{The monodromy representation}

Let $\pi_1(X,x_0)$ be the fundamental group of $X$ with base point
$x_0\in X$. Any local system $\mathcal L$ gives rise to a representation of
$\pi_1(X,x_0)$ on the stalk $\mathcal L_{x_0}$ called the {\it
monodromy representation}. This representation is defined as follows.

Let $\gamma:[0,1]\to X$ be a continuous map,
$\gamma(0)=\gamma(1)=x_0$, representing an element of $\pi_1(X,
x_0)$. By compactness the image of $\gamma$ may be covered with
finitely many open sets $U_1,\dotsc, U_n$ on which $\mathcal L$ is
trivial. And these may be chosen such that $U_{j}\cap
U_{j+1}\ne\emptyset$ and $x_0\in U_n\cap U_1$. Let $x_j\in
U_{j}\cap U_{j+1}$. Then by the definition of locally constant
sheaf we have a sequence of isomorphisms
\begin{equation*}
\mathcal L_{x_0}\overset\sim\leftarrow \mathcal
L(U_1)\overset\sim\rightarrow \mathcal L_{x_1}
\overset\sim\leftarrow \mathcal L(U_2)\overset\sim\rightarrow
\cdots \overset\sim\leftarrow \mathcal
L(U_n)\overset\sim\rightarrow \mathcal L_{x_0},
\end{equation*}
which defines an automorphism (from left to right) of $\mathcal
L_{x_0}$. It follows from the definition of locally constant
sheaf that this automorphism $\mathcal L_{x_0}\to\mathcal
L_{x_0}$ depends only on the homotopy type of $\gamma$. It is
also clear that concatenation of paths corresponds to composition
of maps, giving a representation of $\pi_1(X,x_0)$.

The monodromy representation defines a functor,
 $$
 \left(\begin{array}{cc}\text{local systems on $X$, as}\\
 \text{full subcategory of $\Sh(X)$}\end{array}\right)
 \To \left(\begin{array}{cc}\text{finite dimensional rep's
 of $\pi_1(X,x_0)$,}\\
 \text{and equivariant homomorphisms}\end{array} \right).
 $$
The main observation is that most of the time these two categories
are equivalent (similarly with local systems replaced by locally
constant sheaves and the finite dimensionality constraint on the
right hand side removed).
\begin{prop}
\label{p:monodromyProp} If $X$ has a
universal cover (i.e. $X$ is path connected, locally path
connected and locally simply connected), then the monodromy
functor is an equivalence of categories.
\end{prop}

\begin{proof}
Suppose we are given a representation of $\pi_1(X):=\pi_1(X,x_0)$
on some (finite dimensional) vector space $V$. Let $p:\tilde X\to
X$ be the universal cover of $X$. Then $\pi_1(X,x_0)$ acts on
$\tilde X$ by deck transformations and $p$ is the quotient map.
Define a sheaf by
 $$
 \mathcal L(U):=\left\{\text{locally constant, $\pi_1(X)$-equivariant
 maps $\phi:p\inv(U)\to V$}\right\}
 $$
and the obvious restriction maps (this clearly obeys the sheaf
axiom).

If $x\in X$ and $U$ is a connected, simply connected neighborhood
of $x$, then $p\inv(U)\cong U\x\pi_1(X)$ with $\pi_1(X)$ acting
from the right by right translation. In that case $\mathcal
L(U)\cong V$ by $\phi\mapsto \phi(u,1)$ (independent of $u$ since
$\phi$ is locally constant). The inverse isomorphism comes from
reconstructing $\phi$ from $\phi(u,1)$ using the group action. It
follows that $\mathcal L$ is a locally constant sheaf with stalks
isomorphic to $V$. Moreover, it is easy to see from the
construction of $\mathcal L$ that the monodromy representation is
the representation we started with.
\end{proof}

From now on let us always assume that any topological space $X$
we will consider, if it is connected, has a universal cover.

\begin{rem}[Example]\label{r:exampleMR}  Consider the rank one local
system $\mathcal L^{\alpha}$ on $X=\C^*$ from \eqref{e:Lalpha}.
Here $\pi_1(X)\cong(\Z,+)$ and the monodromy representation
$\rho^\alpha$ of $\mathcal L^\alpha$ is determined by the action
of the generator $1\in \Z$ on a (any) stalk. Then it is easy to
check that $\rho^\alpha(1)=(\exp(2\pi i\alpha))\in GL_1(\C)$,
assuming the generator $1$ is represented by the path $\gamma$ in
$X$ that winds around the origin once in positive, anti-clockwise,
orientation.
\end{rem}

\subsection{Extensions of local systems}
Suppose $U\subset X$ open and $\mathcal L$ is a local system on
$U$. Then one can immediately ask, when is $\mathcal L$ the
restriction of a (unique) local system $\tilde{\mathcal L}$ on
$X$.
 For example for $\mathcal L^{\alpha}\in\Sh(\C^*)$  from
\eqref{e:Lalpha} to come from a local system on
$\C$ it must be constant. This is the case precisely if
$\alpha\in \Z$ as we have already observed.

In general for $\mathcal L$ to extend to $X$ it must have trivial
monodromy around any loop $\gamma$ in $U$ which is contractible
in $X$. Also $U$ should be large enough in $X$ so that
$\pi_1(U)\to\pi_1(X)$ is surjective. Otherwise even if an extension
exists it may not be unique. In a manifold this surjectivity is
automatic if $U$ has complement of codimension $\ge 2$.

\subsection{Outlook~: Perverse sheaves and intersection cohomology}

We have seen in Proposition~\ref{p:monodromyProp} that the
category of local systems on $X$ is completely governed by the
fundamental group $\pi_1(X)$. A common thing one might do given a
local system $\mathcal L$ on $X$ is to compute the cohomology of
$X$ with coefficients in $\mathcal L$. (For example this
corresponds to group cohomology of $\pi_1(X)$ when the universal
cover $\tilde X$ is contractible). When trying to study the
topology of a manifold $X$, the standard cohomology comes from the
constant sheaf. But it is also common to consider the orientation
sheaf to get a more useful theory for non-orientable manifolds.
The idea of choosing a suitable local system adapted to $X$ when
considering cohomology can be generalized to the case when $X$ is
not a manifold but has some singularities.

Suppose for example $X$ is a complex algebraic variety, such as
the nilpotent cone in a Lie algebra. Consider a stratification of
$X$ by smooth subvarieties with one open dense stratum $U$. Then
associated to any local system $\mathcal L$ on $U$ there is a
natural `perverse sheaf' on $X$ extending $\mathcal L$, which
gives rise to an appropriate analogue to cohomology with
coefficients in $\mathcal L$ for the singular space $X$. This
extension in a way assigns to $\mathcal L$ a collection of local
systems in different degrees on the smaller strata by taking into
close account the singularities of $X$ and their severity. For
example, the perverse sheaf extension of $\Const_U$ will not in
general agree with $\Const_X$, unless $X$ is a manifold. And the
cohomology of the perverse extension can vary for two homotopic
but not homeomorphic spaces whereas the cohomology of $\Const_X$
only depends on the homotopy type of $X$. When $\mathcal L$ is
chosen to be the orientation sheaf on $U$ this construction gives
the intersection homology of $X$. For a description of the
perverse sheaf extensions of the
$\mathcal L^\alpha$ from \eqref{e:Lalpha} see
\cite[Section~2.9.14]{KashScha:SheavesBook}.

\section{Lecture - Intersection Homology}
\subsection{Stratified pseudomanifolds}

In this section we will give the definition of intersection
homology of Goresky and MacPherson \cite{GM:IH1}. Some further
references are \cite{Borel:ICBook,Kirwan:IHBook, Borel:ICPaper}.
One wants to deal with topological spaces that are somewhat like
manifolds, but allowed to have singularities. So we will consider
singular spaces that are stratified into a union of manifolds.
\begin{defn}\label{d:StratPM}
An $n$-dimensional {\it stratified pseudomanifold} is a
topological space $X$ with a filtration
 $$
 X=X_0\supset X_1\supset X_2\supset X_3\supset\cdots \supset
 X_n\supset\emptyset
 $$
by closed subsets with the following properties.
\begin{itemize}
\item For any $k$ the stratum $S_k:=X_k\setminus X_{k+1}$ is a topological
manifold of dimension $n-k$, or else empty (the index $k$ may be
thought of as the codimension). In particular $S_n=X_n$ is a
discrete at most countable set.
\item The stratum $S_1=\emptyset$. Therefore $X_1=X_2=:\Sigma$, which we think of
roughly as the singular locus.
\item The open stratum $S_0=X\setminus\Sigma$ is dense in $X$.
\item {\it Local normal triviality along the strata~:} Any $x\in S_k$
has an open neighborhood $U$ and a compact stratified
pseudomanifold of dimension $k-1$, $L=L_0\supset L_2\supset\cdots
\supset L_{k-1}$, called the {\it link} of $x$, with compatible
maps
\begin{equation*}
\begin{array}{rccc}
 & U &\overset\sim\To & \R^{n-k}\x \text{cone}^\circ(L)\\
 &\uparrow & & \uparrow\\
 & U\cap X_i&\overset\sim\To & \R^{n-k}\x \text{cone}^\circ(L_i)\\
 &\uparrow & & \uparrow\\
 & U\cap X_k&\overset\sim\To & \R^{n-k}\x \{\text{vertex}\}\\
\end{array}
\end{equation*}
for any $0\le i<k$. Here $\text{cone}^\circ(L)=L\x [0,\infty
)/(L\x 0)$ is the open cone over $L$ with vertex $L\x 0/(L\x 0)$.
The vertical maps are the obvious inclusions.
\end{itemize}
Note that the indexing of the strata and the $X_k$ by codimension
differs from the one in most references. We use it because it is
more clear when keeping track of dimensions of intersections. An
intersection of a closed submanifold $Z\subset X$ with an $S_k$
is dimensionally transverse if the dimension of $Z$ is reduced by $k$ (or if
it is empty).
\end{defn}

\subsection{Piecewise linear theory}

We begin with the `simplicial' version of intersection homology
theory, which is the historical route and familiar from usual
homology theory. For this we need to consider the category of 
piecewise linear (\pl) topological spaces.

Formally, a {\it {\pl}-space} is a topological space $X$ with a
class $\mathcal T$ of locally finite triangulations such that for
a triangulation $T$ of $X$ in the class any subdivision of $T$
again lies in $\mathcal T$, and also any $T,T'\in \mathcal T$ have
a common refinement in $\mathcal T$. Think for example of a
simplicial complex along with all possible refinements. An
advantage of having a whole class of triangulations is that any
open set $U\subset X$ inherits a \pl-structure. This is useful for
constructing sheaves later on.

A map between \pl-spaces $(X,\mathcal T_X)$ and $(Y,\mathcal
T_Y)$ should be a map $X\to Y$ for which there are triangulations
$T_X\in \mathcal T_X$ and $T_Y\in \mathcal T_Y$ such that the
image of any simplex $\Delta\in T_X$ lies inside a simplex of
$T_Y$.

\begin{defn}\label{d:\pl-pseudomfld} An {\it $n$-dimensional
{\pl}-pseudomanifold} is a \pl-space which is the union of closed
$n$-simplices in some admissible triangulation, and such that any
$(n-1)$-simplex is the face of exactly two $n$-simplices.
\end{defn}

\begin{prop}
Any \pl-pseudomanifold has a filtration making it a stratified
pseudomanifold in the category of \pl-spaces and maps.
\end{prop}

\begin{rem}  The construction of the stratification is straightforward.
Consider the triangulation from
Definition~\ref{d:\pl-pseudomfld}. Then define $X_k$ to be the
union of the closed $(n-k)$-simplices for $k\ge 2$, perforce
$X_1=X_2$. Then $X_k\setminus X_{k-1}$ is obviously a manifold.
That $X\setminus X_2$ is also a manifold follows from the
\pl-pseudomanifold condition that any $(n-1)$-simplex is the face
of exactly two $n$-simplices. We omit the proof of local normal
triviality, but see \cite{Haefliger:BorelEtAl} Proposition 1.4.
\end{rem}

\subsection{Borel-Moore Homology}
Let $X$ be a \pl-space from now on. Borel-Moore homology (also
called homology with closed supports) is to usual homology what
cohomology $H^*(X)$ is to cohomology $H^*_c(X)$ with compact
supports, see for example \cite{Borel:ICBook} 
and \cite[Section 2.6]{GinzChris:book}. 
Following the conventions of \cite{Borel:ICBook} we
will use the otherwise non-standard notation
\begin{equation*}
\begin{array}{rcl} H_*(X) & = &\text{Borel Moore homology of $X$ over $\C$},\\
 H^c_*(X) &= & \text{usual (e.g. singular or simplicial) homology of $X$ over $\C$}.
\end{array}
\end{equation*}
The two notions coincide if $X$ is compact.
\begin{defn}
For  a triangulation $T$ of $X$ an $i$-chain is defined to be a
formal $\C$-linear combination of (possibly infinitely many)
oriented $i$-simplices. Denote by $C_i^T(X)$ the vector space of
all such $i$-chains. For any refinement $T'$ of $T$ there is a map
$C_i^T(X)\to C_i^{T'}(X)$ so that we can take the limit
 $$
 C_i(X)=\lim_{\overset\To {T\in\mathcal T}} C_i^T(X).
 $$
The elements of $C_i(X)$ are called {\it geometric $i$-chains} or
{\it Borel-Moore chains}. The usual boundary map defined on
individual oriented simplices induces a boundary map
$\partial:C_i(X)\to C_{i-1}(X)$.  The {\it Borel-Moore homology}
$H_*(X)$ of $X$ is defined to be the homology of this complex. We
call a geometric $i$-chain $\xi$ with $\partial\xi=0$ an
$i$-cycle.
\end{defn}

A geometric $i$-chain $\xi$ is represented by an $i$-chain in
some $C_i^T(X)$. Its support $|\xi |\subset X$ is the union of
all closed $i$-simplices in $T$ which occur with nonzero
coefficient in $\xi$. It is closed in $X$, which is why
Borel-Moore homology is also often called `homology with closed
supports'.

We may also consider the complex of compactly supported chains
$C_i^c(X)\subset C_i(X)$, which coincides with the limit over
compactly supported chains in $C_i^T(X)$. Since taking homology
commutes with taking (inductive) limits, its homology computes the
usual simplicial homology $H^c_*(X)$ of $X$. The dual complexes
$C^i(X):=C^c_i(X)^*$ and $C_c^i(X):=C_i(X)^*$ then compute the
usual cohomology $H^*(X)$ and its compactly supported version
$H^*_c(X)$, respectively.

\subsection{Poincar\'e duality}
If $X$ is an oriented $n$-manifold, then there 
is the well known Poincar\'e duality pairing 
for the homology of $X$. 
That is, there is a
well-defined intersection product $H_j(X)\x H^c_{k}(X)\to
H^c_{n-j-k}(X)$ such that the resulting map
 $$
H_i(X)\x H^c_{n-i}(X)\to H^c_0(X)\to \C
 $$
gives a perfect pairing. 
Here the final map is the augmentation homomorphism
$H^c_0(X)\to H^c_0(\{\operatorname{pt}\})\cong\C$ which is induced
from mapping $X$ to a point. The following two properties go into
the definition of the intersection product.
\begin{itemize}
\item For any two \pl-chains, it is possible to replace
either one with a homologous one to make its support transversal
to that of the other.
\item
For any two homology classes and transversal chains representing
them, the class of the intersection is independent of the
particular choice of transversal representatives.
\end{itemize}
We can see how Poincar\'e duality fails for pseudomanifolds in an example. 
\begin{rem}[Example]\label{r:cylinders} Let us
consider a cylinder $X=\R\x S^1$ and a
pinched cylinder, $Y=X/(0\x S^1)$, where the circle at the origin
is identified to a single point $p_0$. The homology is easy to
compute~:
\begin{equation*}\begin{array}{lll}
\begin{array}{c|c|c}
i &H_i(X) & H_i^c(X)\\
\hline 0 & 0 & \C\\
1&\C & \C  \\
2&\C & 0
\end{array} &\qquad&
\begin{array}{c|c|c}
i &H_i(Y) & H_i^c(Y)\\
\hline 0 & 0 & \C\\
1&\C & 0  \\
2&\C\oplus\C & 0
\end{array}
\end{array}
\end{equation*}
We see that Poincar\'e duality fails badly for $Y$. Let
$\Sigma:=\{p_0\}\subset Y$. The problem is apparent. Any geometric
$1$-chain representing a generator of $H_1(Y)$ passes through
$p_0$. Hence it cannot be made transversal to $\Sigma$.

Therefore we need to restrict the allowed chains if there is a
singular locus. For any geometric $i$-chain $\xi$ with support
$|\xi |$ we impose the {\it first intersection condition~:}
\begin{equation}\label{e:*p2}
 \dim |\xi |\cap\Sigma \le \dim|\xi |-2.
\end{equation}
For singular surfaces ($2$-dimensional pseudomanifolds) this
condition suffices to define intersection homology and we get~:
\begin{equation*}
\begin{array}{c|c|c}
i &IH_i(Y) & IH_i^c(Y)\\
\hline 0 & 0 & \C\oplus\C\\
1&0 & 0  \\
2&\C\oplus\C & 0
\end{array}
\end{equation*}
\end{rem}

\subsection{Perversities and intersection homology}
So far we have introduced one intersection condition. Supposing
the singular locus $\Sigma$ is as large as it can be for a
pseudomanifold (i.e. codimension $2$), then it just says that 
any allowable chain must intersect $\Sigma$ transversally.  This
however is just the tip of the iceberg, since intersection
homology also takes into account the singularities within the
singular set and so forth. There are different versions of
intersection homology depending on the strictness of the
intersection conditions with the various strata. These are
encoded in what is called the perversity.

\begin{defn}
A {\it perversity} of a fixed dimension $n$ is a map $\p
:\{2,3,\dotsc, n\}\to\Z_{\ge 0}$ such that $ \p(2)=0 $ and
 $$
\p(k)\in\left\{\p(k-1)\, ,\, \p(k-1)+1\right\}.
 $$
The most important perversities are given by
\begin{equation*}
\begin{array}{lcll}
\underline 0&=&(0,0,0,\dotsc\dotsc\dotsc, 0\, ) & \text{$0$-perversity}\\
\underline m&=&(0,0,1,1,2,2,3,\dotsc)\quad & \text{lower middle
perversity}\\
\underline n&=&(0,1,1,2,2,3,3,\dotsc) & \text{upper middle
perversity}\\
\underline t&=&(0,1,2,3,\dotsc, n-2) &\text{top perversity}.
\end{array}
\end{equation*}
Two perversities $\p$ and $\underline q$ are called {\it complementary}
if $\p +\underline q=\underline t$ (e.g. $\underline m$ and
$\underline n$).
\end{defn}

\begin{defn}[Intersection homology]
Let $\p$ be a fixed perversity. A geometric $i$-chain $\xi\in
C_i(X)$ is called $\p$-allowable if
\begin{equation}\label{e:*pk}
\begin{array}{c}
 \begin{array}{lcl}
 \dim(|\xi |\cap X_k)&\le &i - k + \p(k),\\
 \dim(|\partial\xi |\cap X_k)&\le & i-1 -k +\p(k)
 \end{array}\tag*{$(*)_{\p,k}$}
\end{array}
\end{equation}
holds, for all $2\le k\le n$. Notice that for all $\p$,
$(*)_{\p,2}$ coincides with the first intersection condition from
\eqref{e:*p2} for $\xi$ and $\partial\xi$.

Let us denote by
 $$
 I_\p C_i(X)=\{\ \xi\in C_i(X)\ |\ \text{$\xi$ satisfies $(*)_{\p,k}$}\ \}
 $$
the vector space of $\p$-allowable $i$-chains. The boundary map on
geometric chains restricts to $\partial: I_\p C_i(X)\to  I_\p
C_{i-1}(X)$ to define a chain complex. The homology of this
complex is denoted $I_{\p}H_*(X)$ and called the {\it intersection
homology} of $X$ with perversity $\p$.
\end{defn}

\begin{rem}\label{r:PervExamples}
If $\xi\in I_{\underline 0}C_i$, then $|\xi|$ is
transversal to all strata $S_k$ of $X$. This is the strictest of
the intersection conditions. In general if $\p<\underline q$, then  $I_{\p}
C_i(X)\subset I_{\underline q} C_i(X)$ and we have a map
 $$
 I_{\p} H_i(X)\To
I_{\underline q} H_i(X).
 $$

In the weakest case, the top perversity $\underline
t(k)=k-2$, the intersection conditions  $(*)_{\underline t,k}$ on a
chain $\xi$ are
\begin{align*}
 \dim (|\xi |\cap X_k )&\le\dim |\xi | -2, \\
 \dim (|\partial \xi |\cap X_k) &\le\dim |\xi | - 3.
\end{align*}
But for $k>2$ these are automatic consequences of the first
intersection condition, $(*)_{\underline t,2}$. Therefore
$I_{\underline t}C_*(X)$ is just the complex of chains $\xi$ with
$\dim (|\xi |\cap \Sigma)\le \dim |\xi|-2$, and the same condition
for $\partial \xi$.
\end{rem}

\subsection{Intersection homology and normality}
\begin{defn}
An $n$-dimensional stratified pseudomanifold $X$ is called {\it
normal} if any $x\in\Sigma$ has a distinguished neighborhood
$U_x$ such that $U_x\setminus\Sigma$ is connected.
\end{defn}

This definition of normality is a topological analogue of
normality for algebraic varieties. Any pseudomanifold $X$ has a
{\it normalization} $\pi:\tilde X\to X$, where $\tilde X$ is normal,
and $\pi$ is $1$-$1$ on the
open stratum $S_0$ and finite to one on any other stratum.
Explicitly $\#\{\pi\inv(x)\}$ is the number of connected
components of $U_x\setminus\Sigma$. This map $\pi$ also induces a
map on chain complexes
\begin{equation}\label{e:normalizationmap}
 C_*(\tilde X)\to C_*(X).
\end{equation}

\begin{lem} Any $\p$-allowable cycle $\sigma$ of $X$ is the image of a
$\p$-allowable cycle $\tilde \sigma$ of $\tilde X$ under the map
\eqref{e:normalizationmap}.
\end{lem}

\begin{proof}
For $\sigma\in C_i(X)$ there exists a geometric $i$-chain $\tilde\sigma\in C_i(\tilde
X)$ such that $|\tilde\sigma|=\pi\inv(|\sigma|)$. It follows from
$\partial\sigma=0$ that
$|\partial\tilde\sigma|\subset\tilde\Sigma$. So we have
 $$
 |\pi(\partial\tilde\sigma)|\subset |\sigma |\cap\Sigma.
 $$
Here the left hand side is $(i-1)$-dimensional. But intersecting
with $\Sigma$ should cut down the dimension of $\sigma$ by at
least $2$, by the first intersection condition. So
$\partial\tilde\sigma=0$.

Also $(*)_{\p,k} $ for $\sigma$ implies $(*)_{\p,k}$ for
$\tilde\sigma$.
\end{proof}
\begin{cor}\label{c:normal} The map \eqref{e:normalizationmap}
induces an isomorphism
 $$  I_{\p} H_*(\tilde X)\cong I_{\p}H_*(X).
 $$
\end{cor}

The pinched cylinder from Remark~\ref{r:cylinders}, for example,
is not normal. We have that the normalization $\tilde Y$ of the
pinched cylinder consists of two separate hemispheres (the
pinched cylinder is divided into two halves at the singular point
which is split into two). Topologically $\tilde Y$ is a disjoint
union $\R^2\sqcup \R^2$. Thus the result in
Remark~\ref{r:cylinders} giving the intersection homology  of $Y$
is an illustration of Corollary~\ref{c:normal}. From that point
of view, the example of the pinched cylinder just reduces to a
calculation of classical homology. But see Goresky and
MacPherson's Chapter~III in \cite{Borel:ICBook} for a more
sophisticated sample computation.

\subsection{Top- and $\underline 0$-perversity intersection homology}

For the two extremal perversities intersection homology compares
closely with Borel-Moore homology and usual cohomology, respectively. We
compute the cohomology of $X$  by the cochain complex with
$C^i(X)=C^c_i(X)^*$. One can define two maps,
\begin{align}\label{e:chainmap1}
I_{\underline t}C_i(X)&\to  C_i(X),\\
\label{e:chainmap2} C^{n-i}(X)&\to  I_{\underline 0}C_{i}(X).
\end{align}
The first one is just the obvious inclusion. The second one is
defined by Poincar\'e's dual cycle map. Any simplex in the support
of a cochain is replaced by the complementary dimensional
simplices in the barycentric subdivision which meet the center of
the original simplex. See for example
\cite[\S65]{Munkres:AlgTopBook}.  The resulting chain is then
automatically transversal to all the strata, hence lies in
$I_{\underline 0}C_{i}(X)$.  
\begin{prop}\label{p:extremalperv}
Suppose $X$ is oriented and normal. The chain maps
\eqref{e:chainmap1} and \eqref{e:chainmap2} induce isomorphisms
\begin{enumerate}
\item[$(1)$] \ $I_{\underline
t}H_*(X)\overset\sim\to H_*(X)$,
\item [$(2)$]\ $H^{n-i}(X)\overset\sim\to I_{\underline 0}H^*(X)$.
\end{enumerate}
\end{prop}
\begin{rem}
Note that normality implies that $X\setminus X_3$ is still a
topological manifold. The link $L$ of any point in the stratum
$S_2$ must be homeomorphic to $S^1$, and
$\cone^\circ(S^1)\cong\R^2$. So any actual singularities occur in
one codimension higher.

Recall that by Remark~\ref{r:PervExamples}, $\underline
t$-allowability just comes down to the first intersection
condition for $\xi$ and $\partial\xi$. So
Proposition~\ref{p:extremalperv}.$(1)$ basically says that normality
gives the amount of flexibility required to insure any $i$-cycle
$\xi$ is homologous to one transversal enough to $\Sigma=X_2$. See
for example \cite{GM:IH1} for the proof. We will prove $(2)$ in
Remark~\ref{r:0IC}.
\end{rem}

\subsection{Poincar\'e duality pairing}
There is an analogue of the intersection product between intersection homology
groups of different perversities.
\begin{thm}\cite{GM:IH1}\label{t:PD}
Let $X$ be an $n$-dimensional \pl-pseudomanifold.
\begin{enumerate}
\item[$(1)$]
Let $\p,\underline q,\underline r$ be three perversities with
$\p+\underline q\le \underline r$. Then there is a well-defined
pairing
 $$
 \cap~: I_\p H_i(X)\x I_{\underline q} H^c_j(X)\To I_{\underline
 r}H_{n-i-j}^c(X).
 $$
\item[$(2)$]
If $\p+\underline q=\underline t$, then
 $$
  I_\p H_i(X)\x I_{\underline q}
H^c_{n-i}(X)\overset\cap\To I_{\underline
 t}H_{0}^c(X)\To\C.
 $$
defines a perfect pairing.
\end{enumerate}
\end{thm}
\begin{rem}
The intersection pairing is defined by intersecting allowable
cycles that are ``sufficiently transverse''. That is, the
intersection must be in sufficiently general position with respect
to each stratum (to give an element of the correct $I_{\underline
r}C_*$). Then part (1) of the theorem says that any two homology
classes can be represented by cycles with this property (for this
$\underline r$ is not allowed to be too small), and then the
homology class of the resulting intersection is independent of the
choices. If $\p=\underline 0$ and $\underline q=\underline t$, and
$X$ orientable and normal, then Theorem~\ref{t:PD}.(2) is just the
classical cap product $H^{n-i}(X)\x H^c_{n-i}(X)\to \C$ (after
applying Proposition~\ref{p:extremalperv}).

\subsection{Middle intersection homology}\label{s:MiddleIH} If
$X_{k+1}=X_k$ in the filtration of $X$, then the intersection
condition $(*)_{\p,k+1}$ implies $(*)_{\p,k}$. So $I_{\p}C_*(X)$
only depends on values $\p(k)$ for nonempty strata $S_k$. For
example if all strata are of even dimension then the odd values
of $\p$ do not matter. In particular in that case $I_{\underline
m}C_*(X)=I_{\underline n}C_*(X)$, and one can speak of {\it
middle intersection homology}.
\end{rem}

\begin{cor}\label{c:PD}
If $X$ has only even-dimensional strata, then the middle
intersection homology $I_{\underline m}H_*(X)$ satisfies
Poincar\'e duality. In other words
 $$
  I_{\underline m} H_i(X)\x
I_{\underline m} H^c_{n-i}(X)\overset\cap
\To I_{\underline t}H_{0}^c(X)\To\C
 $$
is a non-degenerate pairing.
\end{cor}

\subsection{Sheaves}\label{s:ICsheaves}
One feature of Borel-Moore chains is that for
any open $U\hookrightarrow X$ we have restriction maps
 \begin{align*}
 C_i(X)&\to C_i(U)\\
 I_\p C_i(X)&\to I_\p C_i(U).
 \end{align*}
For compactly supported chains $C_i^c$ one has instead an inclusion
$C_i^c(U)\hookrightarrow C_i^c(X)$.

\begin{defn} Define sheaves  $\dsheaf_X^i$ and $\IpCsheaf_X^i$ on $X$ by
 \begin{align*}
 \dsheaf_X^{-i}(U)&:=C_i(U),\\
 \IpCsheaf_X^{-i}(U)&:= I_\p C_i(U)
 \end{align*}
for all $U\subset X$ open.
These sheaves form cochain complexes $\dsheaf_X^\bullet$ and
$\IpCsheaf_X^\bullet$ of sheaves on $X$, with differential induced
from the boundary maps $\partial$ over each open set.
\end{defn}

The original complexes of vector spaces are recovered (up to the
change in indexing) by applying the functors of global sections or
global sections with compact supports. So we have
\begin{align*}
 H_i(X)&=H^{-i}(\Gamma(\dsheaf_X^\bullet)), &
 H^c_i(X)&=H^{-i}(\Gamma_c(\dsheaf_X^\bullet)), \\
 I_\p H_i(X)&= H^{-i}(\Gamma(\IpCsheaf_X^\bullet)), &
 I_\p H^c_i(X)&= H^{-i}(\Gamma_c(\IpCsheaf_X^\bullet)).
\end{align*}

We also have  sheaves $\mathcal C_X^i$ on $X$ defined by
 $$
 \mathcal C_X^{i}(U)=C^{i}(U):=C^c_{i}(U)^*.
 $$
The corresponding cochain
complex $\mathcal C_X^\bullet$ computes the cohomology of $X$ by
 $$
 H^{i}(X)=H^{i}(\Gamma(\mathcal C_X^\bullet)) \quad
 \text{and} \quad H^{i}_c(X)=H^{i}(\Gamma_c(\mathcal C_X^\bullet)).
 $$
Note that now $\mathcal C_X^i(U)\to C^{T,c}_{i}(U)^*$ for any
admissible triangulation $T$. While the geometric $i$-chains were
a direct limit over all triangulations, the cochains are
necessarily an inverse limit.

\section{Lecture - Derived functors and the derived category of sheaves}

\subsection{Complexes and functors} In the last lecture we introduced
complexes of sheaves $\dsheaf_X^\bullet$ and $\IpCsheaf_X^\bullet$
on a pseudomanifold $X$. These are the objects we will be working
with. Let us for the moment view them as objects in the category
of complexes of sheaves with chain maps. Denote the $i$-th
cohomology (sheaf) of a complex $\mathcal F^\bullet$ by $\mathcal
H^i(\mathcal F^\bullet)$. We are more interested in the cohomology
of our complexes than the complexes themselves.
\begin{defn}
A map of complexes
\begin{equation*}
\begin{array}{ccccccccccc}
\mathcal F^\bullet&~:& \cdots\to& \mathcal F^{m-1}& \to&\mathcal
F^m &
\to& \mathcal F^{m+1}&\to & \mathcal F^{m+2} & \to \cdots \\
\hskip -.2cm q\downarrow& & &\downarrow &    &\downarrow &
&\downarrow&&\downarrow&\\
\mathcal G^\bullet&~:& \cdots \to & \mathcal G^{m-1}&\to &
\mathcal G^m&
\to &\mathcal G^{m+1} &\to &\mathcal G^{m+2}&\to\cdots\\
\end{array}
\end{equation*}
is called a {\it quasi-isomorphism} if it induces isomorphisms
$\mathcal H^i(q)~:\mathcal H^i(\mathcal F^\bullet)\overset\sim\to
\mathcal H^i(\mathcal G^\bullet)$ on the cohomology sheaves of
the complex.
\end{defn}
We would like to work in a category (of complexes of sheaves)
where two quasi-isomorphic complexes
are interchangeable. Moreover this category should
`extend' the category of sheaves in the sense that
that the functors that one uses in the category of
sheaves extend to functors between complexes.

The problem with these two demands is already apparent in the
example of the global sections functor $\Gamma$. Consider the
quasi-isomorphism from Remark~\ref{r:deRham},
\begin{equation*}
\begin{array}{ccccccccccccc}
\Const_X[0]&~:& \cdots\to& 0& \to&\Const_X &
\to&0 &\to &\cdots  & \to &0&\to\cdots\\
\hskip -.2cm q\downarrow& & &\downarrow &    &\downarrow &
&\downarrow&&& &\downarrow&\\
\Omega_X^\bullet&~:& \cdots \to & 0&\to & \Omega^0_X& \to
&\Omega^1_X &\to &\cdots&\to&\Omega_X^{n} &\to  \cdots,
\end{array}
\end{equation*}
and apply $\Gamma$ to it. The cohomology
$H^i(\Gamma(\Omega^\bullet_X))=H^i(X)$ of the lower complex will
usually be nonzero in positive degree, for example if $X$ is a
compact manifold of positive dimension. So $\Gamma(q)$ is far
from being a quasi-isomorphism. However we do recover the identity
$\Gamma(\Const_X)=H^0(X)$ in degree zero.

These properties reflect the fact that $\Gamma$ is left exact and
not exact. The correct definition of the extension of a left
exact functor to complexes may be illustrated by the above
example. It says not to apply $\Gamma$ to arbitrary complexes,
but only to ones whose objects are `$\Gamma$-acyclic'. The
$\Omega_X^i$ are examples of $\Gamma$-acyclic sheaves. They are
soft for example (see Section~\ref{s:acyclics}). $\Gamma$-acyclic
sheaves have good local to global behaviour, making $\Gamma$ act
like an exact functor.\footnote{This may sound contradictory,
since for example $\Omega_X^{n-1}\overset{d_n}\To\Omega_X^{n}\to
0$ is exact for an $n$-manifold $X$, while the global sections
can give cohomology in degree $n$. But that can happen because
the kernel of $d_n$ need not be $\Gamma$-acyclic.} Any other
complexes such as $\Const_X[0]$ should be replaced by
quasi-isomorphic ones such as $\Omega_X^\bullet$ before applying
$\Gamma$. Let us now explain these definitions more formally.

\subsection{Right derived functors and $\mathcal D^+(X)$}
Recall that an  object $I$ in an abelian category is called {\it
injective} if $Hom(\ \ , I)$ is an exact functor.
Injective objects have the following useful property.
\begin{itemize}
\item
If $q:I^\bullet\to J^\bullet$ is a quasi-isomorphism between
complexes of injectives, then there is a quasi-isomorphism
$p:J^\bullet\to I^\bullet$ such that $p\circ q:I^\bullet\to
I^\bullet$ is homotopic to the identity.
\end{itemize}
This implies that any additive functor $F$ will act on injectives as if
it were an exact functor, that is, preserving quasi-isomorphisms.

Let us introduce a slightly better category to work in, the
homotopy category
\begin{equation*}
\mathcal K^+(\Sh(X))=\left( \begin{array}{c} \text{bounded below cochain
complexes of sheaves in $\Sh(X)$}\\
\text{with cochain maps modulo homotopy equivalence}
\end{array}\right ).
 \end{equation*}
Any two quasi-isomorphic complexes of injectives are isomorphic
in $\mathcal K^+(\Sh(X))$. The (bounded below) derived category
of $\Sh(X)$ is defined by
 $$
\mathcal D^+(X):=\mathcal K^+(\Sh(X))[\ (\text{quasi-isomorphisms})\inv\ ],
 $$
adjoining an inverse to any quasi-isomorphism in the homotopy
category (localization). The category $\mathcal D^+(X)$ is no
longer abelian but instead a `triangulated category'. See for
example \cite{Iversen:SheafCohBook} or \cite{Weibel:HomAlgBook}.

Let $\mathcal K^+(Inj(X))\subset\mathcal K^+(\Sh(X))$ be the full
subcategory whose objects are complexes of injectives. Then we
also have an inclusion $\iota:\mathcal
K^+(Inj(X))\hookrightarrow\mathcal D^+(X)$. To formally define a
right derived functor of a functor $F$, we need a functor
$I:\mathcal D^+(X)\to\mathcal K^+(Inj(X))$ assigning to any
complex an injective resolution,
\begin{equation}\label{e:I}
\begin{array}{rccccccccc}
\mathcal F^\bullet~:& \cdots\to& 0& \to&\mathcal F^0
 &\to& \mathcal F^1&\to & \mathcal F^2 & \to \cdots\ \\
\downarrow \quad\ & &\downarrow &    &\downarrow &  &\downarrow&&\downarrow&\\
I(\mathcal F^\bullet)~:& \cdots \to & 0&\to & \mathcal I^0&\to
&\mathcal I^1 &\to &\mathcal I^2&\to\cdots,
\end{array}
\end{equation}
such that $I$ and $\iota$ are inverse equivalences of categories.
Such a functor exists for the category of sheaves, as we will see
in the next section.

\begin{defn}
Suppose that $F:\Sh(X)\to\Sh(Y)$ is a left exact functor. Then
define
$RF(\mathcal F^\bullet):=F(I(\mathcal F^\bullet))$
for a complex of sheaves $\mathcal F^\bullet$.
This defines a functor
$$
RF:\mathcal D^+(X)\to\mathcal D^+(Y)
$$
called a {\it right derived functor} of $F$.
The $i$-th right derived functor
$$
R^i F:\mathcal D^+(X)\to\mathcal \Sh(Y)
$$
is defined by $R^i F=\mathcal H^i\circ RF$.
\end{defn}

Given two different injective resolutions of a complex
one can always construct a quasi-isomorphism between them and
uniquely up to homotopy, rendering them canonically isomorphic in
the homotopy category. So the choice of the functor $I$ hardly
matters.

Suppose $\mathcal F[0]$ is the complex with $\mathcal F$ in
degree $0$ and zero elsewhere. If $F$ is a  left exact functor,
then applying $F$ to the quasi-isomorphism \eqref{e:I} induces an
isomorphism $F(\mathcal F)\cong R^0 F(\mathcal F[0])$. This is
the sense in which $RF$ is an extension of the functor $F$ to
complexes.

\subsection{Injectives in  $\Sh(X)$.}
In the category of vector spaces, that is the case of
$\Sh(\{pt\})$, every object is injective and projective (the dual property)
by an application of Zorn's Lemma.
More generally there is the following
characterization of injective sheaves due to N.~Spaltenstein.
\begin{lem}
\label{l:Spaltenstein}
A sheaf $\mathcal F\in\Sh(X)$ is injective precisely if for all
open subsets $V\subset U$, the restriction map
$r_{UV}:\mathcal F(U)\to\mathcal F(V)$ is a surjection.
\end{lem}
This says that for sheaves of vector spaces being injective
is equivalent to being {\it flabby} (any section
comes from a global section).
See \cite[V~Lemma~1.13]{Borel:ICBook} for a more
general formulation and the proof of this lemma.
In the setting of
sheaves of modules over a ring, say, one needs also the further
requirements that  the
stalks and sections of $\mathcal F$ should be injective
and the restriction maps be split.

The category of sheaves has ``enough injectives'', meaning that
every sheaf embeds into an injective one. Explicitly, we may let
$\mathcal F\hookrightarrow \mathcal I$ where
 $$
 \mathcal I(U):=\prod_{x\in U} \mathcal F_x
 $$
with the obvious restriction (by projection) maps. Embedding the
cokernel of this inclusion into another injective and so forth,
one obtains an injective resolution of $\mathcal F$ called the
{\it Godement resolution}. So any sheaf has an injective
resolution. Viewed inside the derived category this resolution
becomes an isomorphism (i.e. quasi-isomorphism of complexes)
\begin{equation*}
\begin{array}{ccccccc}
 \cdots\to& 0& \to&\mathcal F &\to& 0&\to\cdots \\
 &\downarrow &    &\downarrow &  &\downarrow&\\
 \cdots \to & 0&\to & \mathcal I^0&\to &\mathcal I^1 &\to \cdots\\
\end{array}
\end{equation*}
The existence of enough injectives also implies that any complex
of sheaves $\mathcal F^\bullet$ is quasi-isomorphic to a complex
of injectives $I(\mathcal F^\bullet)$. This quasi-isomorphic
complex of injectives may be obtained as the total complex of a
double complex constructed from the injective resolutions of all
of the individual $\mathcal F^i$'s (the {\it Cartan-Eilenberg}
resolution).

\subsection{The bounded derived category}
For our purposes it will be more convenient to consider complexes
that are bounded in both directions. These form a full
subcategory $ \mathcal D^b(X) $ in $\mathcal D^+(X)$ called the
{\it bounded derived category} of sheaves on $X$. One constructs
bounded complexes from unbounded ones by {\it truncation}: Let
$\tau_{\le k}~:\mathcal D^+(X)\to \mathcal D^b(X)$ be defined on
objects by
 $$
 \tau_{\le k}(\mathcal F^\bullet)= (\cdots \to \mathcal
 F^{k-2}\to\mathcal F^{k-1}\to\ker(d)\overset d\to 0\to
 0\to\cdots)
 $$
with the obvious extension to morphisms. $\tau_{\le k}$ is called
the (right) truncation functor at $k$ since
\begin{equation*}
\mathcal H^i(\tau_{\le k}(\mathcal F^\bullet))=\begin{cases}
 \mathcal H^i(\mathcal F^\bullet) & i\le k,\\
 0 & i>k.
 \end{cases}
\end{equation*}

To define derived functors between bounded derived categories we
need bounded injective resolutions. These may be obtained as
follows. If $\mathcal I^\bullet \in\mathcal D^+(X)$ is an
injective resolution of $\mathcal F^\bullet\in \mathcal D^b(X)$
and $X$ has dimension $n$, then $ \tau_{\le n+1}\mathcal
I^\bullet$ is still an injective resolution of $\mathcal
F^\bullet$, see \cite[V Proposition 1.17]{Borel:ICBook}.

Let us now introduce some explicit functors 
between derived categories of sheaves.

\subsection{Push-forwards} Let $f:X\to Y$ be a continuous
map. Most interesting functors on sheaves are
left exact~:
\begin{align*}
\Gamma,\Gamma_c:&\Sh(X)\to \Vect, \\
f_*, f_!~:& \Sh(X)\to\Sh(Y).
\end{align*}
Here the first two functors are the global sections functors
introduced in Section~\ref{s:restrictions}. The second are the
push-forward (or direct image) functor $f_*$ and the push-forward
with proper supports $f_!$. They are defined by
\begin{align*}
 f_*\mathcal X (U)&=\mathcal X(f\inv(U)),\\
 f_!\mathcal X(U)&=\{s\in\mathcal X(f\inv(U)) \ |\ f:\Supp(s)\to U
 \text{\ is proper}\},
\end{align*}
where $\mathcal X$ is a sheaf on $X$.
Recall that a map between topological spaces is called proper if
the inverse image of a compact set is compact. For a composition
of maps $f\circ g$ the push-forwards are $(f\circ g)_*=f_*\circ
g_*$  and $(f\circ g)_!=f_!\circ g_!$ respectively. If
$f=\pi:X\to\{pt\}$, then we recover $ \pi_*=\Gamma$ and
$\pi_!=\Gamma_c$. We have the derived functors
 $$
 Rf_*, Rf_!:\mathcal D^b(X)\to \mathcal D^b(Y).
 $$

In the special case where $f=j:X\hookrightarrow Y$ is an
inclusion of a locally closed subset,
 $$
(\open_!\mathcal X)_y=\begin{cases} 0 & y\notin X\\ \mathcal X_y&
y\in X\end{cases}.
 $$
This functor is therefore called {\it extension by zero}, and it
is exact. The extension by zero of the constant sheaf
$\Const_{\{x_0\}}$ from $\{x_0\}\hookrightarrow X$ is the
skyscraper sheaf $\mathcal S$ from Remark~\ref{r:stalks}.

\subsection{Pull-backs}\label{s:pullbacks}
Any continuous map $f:X\to Y$ gives rise to a
pull-back (or inverse image) functor
 $
 f^*:\Sh(Y)\to \Sh(X)
 $,
which takes any $\mathcal Y\in \Sh(Y)$ to the sheaf on $X$
associated to the presheaf
 $$
U\mapsto \lim_{\overset{\longrightarrow}{\{V\,\text{open}\,
|\,f(U)\subset V\}} }\mathcal Y(V).
 $$
The key property of the pull-back is that
 $$
 (f^*\mathcal Y)_x=\mathcal Y_{f(x)}.
 $$
The restriction functors from Section~\ref{s:restrictions} are
pull-backs in the special case where $f$ is an inclusion of a
locally closed subset. The pull-back $f^*$
defines a functor
 $$
 f^*:\dsheaf^b (Y)\to\dsheaf^b(X)
 $$
directly without needing to be derived, since it is exact.
 For a composition of maps $f\circ g$ the
pull-back is $(f\circ g)^*=g^*\circ f^*$.

Pull-back and push-forward are adjoint functors. I.e. 
there is a canonical isomorphism
 $$
\Hom_{\Sh(X)}(f^* \mathcal Y,\mathcal
X)\cong\Hom_{\Sh(Y)}(\mathcal Y,f_*\mathcal X).
 $$
For example the adjunction morphism $1\to f_*\circ f^*$ is given
explicitly (on $\mathcal Y$) by the natural maps $\mathcal Y(V)\to
f^*(\mathcal Y)(f\inv(V))$. In the derived category $f^*$ becomes
left adjoint to $Rf_*$ using the above adjunction morphism
together with the natural transformation coming from \eqref{e:I}.

\subsection{Shifts}
Shifting the indices of a complex to the right defines a functor
$$
[n]:\mathcal D^b(X)\to \mathcal D^b(X).
$$
This functor takes
$\mathcal F^\bullet\mapsto \mathcal F^\bullet[n]$,
where $(\mathcal F^\bullet[n])^i=\mathcal F^{i+n}$ and the
differential of the complex is
conventionally multiplied by a factor of $(-1)^n$.

\subsection{Computing derived functors}\label{s:acyclics}
Injective sheaves are a good theoretical tool, but not always
practically useful.
For computing right derived functors it tends to be more convenient to use
resolutions that are not quite injective but still in a sense `flexible enough'
for the functor in question.

\begin{defn} Let $F$ be a left exact functor from $\Sh(X)$ to $\Sh(Y)$. 
A sheaf $\mathcal A\in\Sh(X)$ is called {\it $F$-acyclic} if
$R^iF(\mathcal A)=0$ for all $i\ne 0$.
\end{defn}

By the spectral sequence
\begin{equation*}
E^{ij}_{2}= H^i(R^j F(\mathcal A^\bullet))\Rightarrow
R^{i+j}F(\mathcal A^\bullet)
\end{equation*}
for computing the cohomology of the Cartan-Eilenberg double
complex we see that $F$-acyclic complexes are good enough for
computing $RF$ (up to isomorphism). In other words, if $\mathcal
F^\bullet$ is quasi-isomorphic to an $F$-acyclic complex $\mathcal
A^\bullet$, then $RF(\mathcal F^\bullet)\cong F(\mathcal
A^\bullet)$ in $\dsheaf^b(X)$.

\begin{defn} A sheaf $\mathcal F$ is called {\it soft}
if the restriction maps
 $$
 \Gamma(U,\mathcal F)\to\Gamma(K,\mathcal F)
 $$
are surjective for all compact
$K$ and open $U\subset X$ containing $K$. Here $\Gamma(K,\mathcal
F):=\Gamma(\closed^*\mathcal F)$ for the inclusion
$\closed:K\hookrightarrow U$.
\end{defn}

In some references the definition of soft requires the
surjectivity for all closed $K$, in which case the above version
is called c-soft. However if $X$ is a union of countably many
compact sets (which is always the case for us) then the two
definitions are equivalent.

A soft sheaf $\mathcal F\in \Sh(X)$ is automatically
$\Gamma_c$-acyclic and $\Gamma$-acyclic, and the restriction of a
soft sheaf to a locally closed subset is again soft. See for example
\cite[III~Theorem~2.7 and IV~Theorem~2.2]{Iversen:SheafCohBook} or
\cite[II 2.5]{KashScha:SheavesBook}.

\begin{lem}$\IpCsheaf_X^\bullet$ (and $\mathcal C_X^\bullet,
\mathcal D_X^\bullet$) are complexes of soft sheaves.
\end{lem}

\begin{proof}
Let $K$ and $U\subset X$ be as in the definition of softness
above. Suppose  $s\in\Gamma(K,\IpCsheaf_X^{-i})$ is represented by
$\xi\in\Gamma(V,\IpCsheaf_X^{-i})$ for some open neighborhood $V$
of $K$. We may assume $V\subset U$. The problem with $\xi$ is that
it might not have closed support in $U$. Choose a triangulation of
$V$ fine enough such that there is a closed \pl-neighborhood $N$
of $K$ entirely in $V$. Then the intersection of  $\xi$ with $N$
lies in $ \Gamma(U,\IpCsheaf_X^{-i})$ and also represents $s$, and
we are done. The proofs for $\mathcal C_X^\bullet$ and  $\mathcal
D_X^\bullet$ are similar.
\end{proof}

As a consequence we see that the intersection homology is the {\it
hypercohomology} (derived functor of $\Gamma$) of the complex
$\IpCsheaf_X^\bullet$,
\begin{align*}
I_{\p}H_i(X)&=R^{-i}\Gamma(\IpCsheaf^\bullet_X)=:\mathbb
H^{-i}(X,\IpCsheaf_X^\bullet),\\
I_{\p}H^c_i(X)&=R^{-i}\Gamma_c(\IpCsheaf^\bullet_X)=:\mathbb
H_c^{-i}(X,\IpCsheaf_X^\bullet).
\end{align*}

The complex of cochains  $\mathcal C_X^\bullet$ on a
\pl-pseudomanifold $X$ introduced in the end of
Section~\ref{s:ICsheaves} is a soft resolution of the constant
sheaf $\Const_X[0]$. Therefore the cohomology of $\Gamma(\mathcal
C_X^\bullet)$ is the sheaf cohomology (hypercohomology of a single
sheaf) of $\Const_X$,
\begin{equation*}
H^i(X)=R^{i}\Gamma(\mathcal C^\bullet_X)=\mathbb
H^i(X,\Const_X[0])=: H^{i}(X,\Const_X).
\end{equation*}

\vskip.2cm There are many useful formulas and identities for
dealing with derived functors. We mention below some of the ones
that will be needed later on, see for example V~\S 10 in
\cite{Borel:ICBook}.
\subsection{Push-forward with proper supports} For the projection $\pi:X\to\{\pt\}$
and  a complex of sheaves $\mathcal X^\bullet$ is on $X$ we have
$R^i\pi_!(\mathcal X^\bullet)=\mathbb H_c^i(X,\mathcal
X^\bullet)$, by definition. Now let $f:X\to Y$ be any continuous
map and $\mathcal X$ a sheaf on $X$. Then there is a natural
isomorphism
 $$
 (f_!\mathcal X)_y\cong \Gamma_c(f\inv(y),\mathcal X|_{f\inv(y)})
 $$
for the stalk of the push-forward with proper supports. Applying
this isomorphism for the sheaves in $\mathcal I^\bullet=I(\mathcal
X^\bullet)$ in the definition of $Rf_!(\mathcal X^\bullet)$ gives
the formula
 \begin{equation}\label{e:ProperStalks}
(R^i f_!\mathcal X^\bullet)_y\cong H^i(\Gamma_c(f\inv(y),\mathcal
I^\bullet|_{f\inv(y)})) \cong \mathbb H^i_c(f\inv(y),\mathcal
X^\bullet|_{{f\inv(y)}}).
 \end{equation}
Here the second isomorphism holds since $\mathcal
I^\bullet|_{f\inv(y)}$ is  a soft resolution of $\mathcal
X^\bullet|_{f\inv(y)}$.

\subsection{Composition of derived functors}
The push-forward functor $f_*$ takes injectives to injectives (for
formal reasons, since it is right adjoint to the exact functor
$f^*$). Therefore $R{(g\circ f)_*}\cong Rg_*\circ Rf_*$ for maps
$f:X\to Y$ and $g:Y\to Z$. More generally, $R(G\circ F)\cong
RG\circ RF$ if $F$ turns injectives into $G$-acyclic objects.

\subsection{Base change}\label{s:BaseChange}
Suppose we have a Cartesian square, that is, a commutative diagram
\begin{equation*}
\begin{array}{rcl}
X_1&\overset {h_1}\To &X_2\\
f_1\downarrow\  & &\ \downarrow f_2\\
Y_1&\overset {h_2}\To&Y_2\quad
\end{array}
\end{equation*}
such that $(f_1,h_1):X_1\overset\sim\to\{(y,x)\in Y_1\x X_2\ |\
h_2(y)=f_2(x)\}$. Then there is a natural isomorphism of functors
$\mathcal D^b(X_2)\to\mathcal D^b(Y_1)$,
 $$
 h_2^*\circ R(f_2)_!= R(f_1)_!\circ h_1^*.
 $$

\section{Lecture - Local intersection cohomology and Deligne's construction}

\subsection{Intersection cohomology}
We will denote the hypercohomology of the complex
$\IpCsheaf^\bullet_X$ from Section~\ref{s:ICsheaves} by
 $$
 I_\p H^k(X):=\mathbb H^k(X,\IpCsheaf_X^\bullet),
 $$
and call it {\it intersection cohomology}. This differs from
intersection homology only in the indexing: $I_\p H^k(X)=I_\p
H_{-k}(X)$. Note that in these conventions intersection cohomology
is concentrated in negative degrees. In the next few sections we
want to investigate intersection cohomology locally in the
neighborhood of a point.
\begin{defn}
Let $\mathcal F^\bullet\in\dsheaf^+(X)$ and $x\in X$. Then the
{\it stalk cohomology} functor at $x$ is the functor $ \mathcal
H^i_x~:\dsheaf^+(X)\to \Vect$ defined by
\begin{equation*}
\mathcal H^i_x(\mathcal F^\bullet):=H^i(\mathcal F^\bullet_x).
\end{equation*}
An equivalent description is
 \begin{equation}\label{e:stalks}
\mathcal H^i_x(\mathcal F^\bullet)=\mathcal H^i(\mathcal
F^\bullet)_x=\lim_{\overset{\longrightarrow}{\{\,\text{$U$ open}\,
|\,\ x\in U\}}}H^i(\mathcal F^\bullet(U)).
 \end{equation}
\end{defn}

\begin{rem} Intersection cohomology, unlike usual cohomology,
{\it can} distinguish between spaces that are homotopic 
(but not homeomorphic) to one another. For example, as we shall see, 
even though every point in $X$ has a contractible neighborhood,
the local intersection cohomology and the stalks
 \begin{equation}\label{e:localIC}
 \mathcal H^{-i}_x(\IpCsheaf_X^\bullet)=
\lim_{\overset{\longrightarrow} {\{U\ \text{open}\ |\  x\in U\}}}
I_\p H_i(U)
 \end{equation}
need not be trivial. Accordingly,
$\IpCsheaf^\bullet_X$ is also not generally quasi-isomorphic to a
single sheaf (unlike $\mathcal C^\bullet_X$).
\end{rem}

Recall the notations  from Definition~\ref{d:StratPM} for the
stratified pseudomanifold structure on $X$.  Let $x$ lie in a
stratum $S_k$. Then we know that $x$ has a distinguished
neighborhood $U\cong \R^{n-k}\x \cone^\circ(L)$, where $L$ is a
$(k-1)$-dimensional compact pseudomanifold. To compute the stalk
cohomology \eqref{e:localIC}  it suffices to consider such
distinguished neighborhoods $U$, since they form a neighborhood basis.
Now $I_\p H_i(U)$ can be computed in terms of $L$ in two steps.

\subsection{Suspension}
If $X$ is a $n$-dimensional stratified pseudomanifold, then $\R\x
X$ is an $(n+1)$-dimensional stratified pseudomanifold with
respect to the filtration $(\R\x X)_i=\R\x X_i$. There is a
natural suspension map $$ \R \x\underbar{\ }~:C_i(X)\to
C_{i+1}(\R\x X) $$ on geometric chains satisfying $|\R\x\xi|=\R\x
|\xi |$ on supports. Suppose $\p$ is an $(n+1)$-dimensional
perversity giving rise to intersection conditions on both $X$ and
$\R\x X$. It is easy to see that the suspension $\R\x \xi$ of a
$\p$-allowable chain is again $\p$-allowable. We have the
following `K\"unneth formula', see \cite[Chapter II \S2]{Borel:ICBook}.

\begin{prop}\label{p:suspension} The suspension map
$$\R \x\underbar{\ }~:C_i(X)\to C_{i+1}(\R\x X)$$
induces an isomorphism $I_\p H_i(X)\overset \sim\to I_\p H_{i+1}(\R\x X)$.
\end{prop}

\subsection{Intersection homology of a cone}
Let $L\supset L_2\supset\cdots\supset L_{k-1}$ be a compact
$(k-1)$-dimensional stratified pseudomanifold, and $Y=\cone^\circ
(L)$. So $Y$ is a stratified pseudomanifold with
\begin{equation*}
 Y_i=\begin{cases}\cone^\circ(L_i) & 2\le i<k,\\
 \{\, v\, \}& i=k,
\end{cases}
\end{equation*}
where $v$ is the vertex of the cone. There is a natural map on
chains, $$\cone^\circ(\underbar{\ }\,):C_i(L)\to
C_{i+1}(\cone^\circ(L)),$$ satisfying $|
\cone^\circ(\xi)|=\cone^\circ(|\xi |)$ on supports. Let $\p$ be a
fixed perversity of dimension $k$.
\begin{lem}
If $\xi\in I_\p C_{i-1}(L)$, then
 $$
 \cone^\circ(\xi) \in I_\p C_i(\cone^\circ(L))\  \iff \begin{cases}
 \text{either }& i>k-\p(k),\\
 \text{or }& i=k-\p(k) \text{ and }\partial\xi=0.
 \end{cases}
 $$
\end{lem}
\begin{proof} Note that $\cone^\circ(\xi)$
always meets the vertex, in other words
the most singular stratum of $Y=\cone^\circ(L)$.
So the dimension $i$ of $\cone^\circ(\xi)$ must be big enough for this to be
$\p$-allowable. That is,
\begin{align}
0&\le \dim (|\cone^\circ(\xi)|)-k+\p(k),\\
 0&\le\dim(|\partial(\cone^\circ(\xi))|)-k+\p(k).
\end{align}
This implies the necessary direction of the lemma. All other $\p$-allowability
restrictions for $\cone^\circ(\xi)$ follow from the ones for $\xi$.
\end{proof}

\begin{prop}\label{p:coning} The map
$\cone^\circ(\underbar{\ }\,):C_i(L)\to C_{i+1}(\cone^\circ(L))$ from above
induces a quasi-isomorphism
 $$I_\p
C_\bullet(L)\overset \sim\To \tau_{\le{\p(k)}}
(I_\p C_{\bullet+1}(\cone^\circ(L))).
 $$
\end{prop}
For a proof of this proposition see \cite[Chapter II \S3]{Borel:ICBook}.

\subsection{Local intersection cohomology}
Let $x\in S_k\subset X$ and $U\cong \R^{n-k}\x \cone^\circ(L)$ be
a distinguished neighborhood. Combining
Propositions~\ref{p:suspension} and \ref{p:coning} one has
\begin{equation*}
I_\p H^i(U)= \begin{cases}I_\p H^{n-k+1-i}(L) & -i\le -n+\p(k),\\
 0 & -i>-n+\p(k).\end{cases}
\end{equation*}
If $V\subset U$ is a compatible distinguished neighborhood of $x$, then
the restriction map
 $$
\mathcal I_\p \mathcal C^\bullet_X (U)\to\mathcal I_\p\mathcal
C^\bullet_X (V)
 $$
is a quasi-isomorphism, essentially since all cohomologically
nontrivial cycles come from $L$. Therefore
 $$
\IpCsheaf_X^\bullet (U)\to(\IpCsheaf^\bullet_{X})_x.
 $$
is a quasi-isomorphism and we have computed the stalk cohomology
of $ \mathcal I_\p \mathcal C^\bullet_X $.
\begin{cor} Let $x\in X$ be a point lying in a stratum $S_k$ with link $L$.
The stalk cohomology $\mathcal H^i_x(\mathcal I_\p\mathcal
C_X^\bullet)$ is concentrated in degrees $-n\le i\le -n+\p(k)$,
where it is determined by the intersection homology of $L$ as
shown in the table below.
\begin{equation*}
\begin{array}{c||c|c|c|c}
i &-n & -n+1 & \cdots & -n+\p(k)
\\
\hline \mathcal H^i_x(\mathcal I_\p\mathcal C_X^\bullet) & I_\p
H_{k-1}(L) & I_\p H_{k-2}(L) &\quad \cdots\quad &I_\p
H_{k-\p(k)-1}(L)
\end{array}
\end{equation*}
\end{cor}

\subsection{A new beginning} \label{s:towardsDeligne}
The rest of this lecture will be taken up with a recursive
stratum by stratum  construction of a new complex of sheaves
$\mathcal P^\bullet$ with which to compute intersection
cohomology. This construction is due to Deligne and appears in
\cite{GM:IH2}. An immediate advantage the definition of $\mathcal
P^\bullet$ has over $\IpCsheaf^\bullet_X$ is that it does not
depend on a \pl-structure on $X$.

Let $X$ be a fixed $n$-dimensional stratified pseudomanifold with
all the usual notations. So $X\supset X_2\supset
X_3\supset\cdots\supset X_n$, $\Sigma=X_2$ is the `singular
locus' and $S_k=X_k\setminus X_{k-1}$ is the codimension $k$
stratum. Also define open sets $$ U_k:=X\setminus X_k, \qquad
2\le l\le n. $$ The smallest of these, $U_2$, agrees with the
open stratum $S_0$. We have inclusions
 $$
S_k\overset{\closed_k}\hookrightarrow
U_{k+1}\overset{\open_k}\hookleftarrow U_k,\qquad\text {for $2\le
k<n$}
 $$
where $\closed_k$ is a closed embedding and $\open_k$ an open
embedding.

\subsection{Extensions of complexes and the attachment map}
Let us consider the case of an open embedding, $\open
:V\hookrightarrow W$, later to be taken to be $\open_k$ from
above. Suppose $\mathcal F_V^\bullet\in\mathcal D^b(V)$. Then one
possible extension of $\mathcal F_V^\bullet$ to $W$ is given by
$R\open_*(\mathcal F_V^\bullet)$. Here $R\open_*(\mathcal
F_V^\bullet)$ is an extension in the sense that $\open^*
(R\open_*(\mathcal F_V^\bullet))\cong \mathcal F_V^\bullet$.

If $\mathcal F_W^\bullet$ is any other extension of $\mathcal
F_V^\bullet$, then it maps to this one by what is called the {\it
attachment map},
 $$
\alpha:\mathcal F_W^\bullet\to R\open_*(\mathcal F_V^\bullet).
 $$
This is just the adjunction morphism $1\to R\open_*\circ \open^*$
from Section~\ref{s:pullbacks} applied to $\mathcal F_W^\bullet$.
Explicitly, for any $U\subset W$ open, $R\open_*(\mathcal
F_V^\bullet)(U)=I(\mathcal F_V^\bullet)(U\cap V)$, and $\alpha(U)$
is the composition of the natural map $\mathcal F_W^\bullet(U)\to
I(\mathcal F_W^\bullet)(U)$ with the restriction map from $U$ to
$U\cap V$, followed by the isomorphism $\mathcal
F_W^\bullet|_V\to\mathcal F_V^\bullet$. On $V$ the attachment map
$\alpha$ is clearly an isomorphism.

\subsection{The Deligne sheaf}
The restriction $\mathcal I_\p\mathcal C_X^\bullet|_{U_2}$ to the
open stratum of $X$ is always quasi-isomorphic to the (shifted)
orientation sheaf $\mathbb O_{U_2}[n]$, where $n$ is the
dimension of $X$. This is simply because there are no
intersection conditions on $U_2$, and we are in the case of usual
Borel-Moore chains. And since the isomorphism $H^n(\R^n)\cong \C$
in Borel-Moore homology depends on a choice of orientation of
$\R^n$, it follows that the local Borel-Moore homology $\mathcal
H^{-n}(\mathcal C_{U_2}^\bullet)$ on the manifold $U_2$ is the
orientation sheaf. The other local cohomology groups vanish of
course, as they do for $\R^n$.

\begin{defn}
Let $\mathcal L$ be a local system on the open stratum $U_2\subset X$ and $\p$ a
fixed perversity.
Define $\mathcal P_{k}^\bullet=\mathcal P_k^\bullet(\p,\mathcal L)\in \mathcal D^b(U_k)$
inductively as follows.
\begin{itemize}
\item [(1)]
Set $\mathcal P_{2}^\bullet:=\mathcal L[n]$.
\item [(2)]
For $k\ge 2$ set $ \mathcal P_{k+1}^\bullet:=\tau_{\le
-n+\p(k)}R{\open_k}_* (\mathcal P_{k}^\bullet). $
\end{itemize}
The complex $\mathcal P_{\p}^\bullet(\mathcal L):=\mathcal
P^\bullet_{n+1}$ is called the {\it Deligne sheaf} corresponding
to $\mathcal L$ and $\p$. If $\mathcal L$ is the orientation
sheaf $\mathbb O_{U_2}$ we may also write $\mathcal P^\bullet_\p$
for $\mathcal P^\bullet_{\p}(\mathbb O_{U_2})$.
\end{defn}

\begin{prop}
The attachment map corresponding to $\open_k:U_k\to U_{k+1}$ gives
rise to a quasi-isomorphism
 $$
\alpha_k:\mathcal I_\p\mathcal C^\bullet|_{U_{k+1}}\To \tau_{\le
-n+\p(k)} R{\open_k}_*(\mathcal I_\p\mathcal C^\bullet|_{U_k}). $$
\end{prop}

\begin{proof}
Since $\mathcal I_\p\mathcal C^\bullet|_{U_k}$ is soft we have
that
 $$
{\open_k}_*(\mathcal I_\p\mathcal C^\bullet|_{U_k})\to
R{\open_k}_*(\mathcal I_\p\mathcal C^\bullet|_{U_k})
 $$
is a quasi-isomorphism. Recall that $\open_k:U_k\hookrightarrow
U_{k+1}=U_k\sqcup S_k$. Let $x\in S_k$ and let $V\subset U_{k+1}$
be a distinguished neighborhood with
\begin{equation*}
\begin{array}{ccc}
V &\overset \sim\To &\R^{n-k}\x\cone^\circ(L)\\
\uparrow & & \uparrow\\
V\cap U_k &\overset\sim\To &\R^{n-k}\x (\cone^\circ(L)\setminus\{v\}).
\end{array}
\end{equation*}
Note that topologically $\cone^\circ(L)\setminus\{v\}\cong \R\x L$.
By a result analogous to Proposition~\ref{p:suspension}
(see Chapter~II, Proposition~3.4 in \cite{Borel:ICBook})
we have isomorphisms
$$
I_\p H_i(L)\to I_\p H_{i+1}(\cone^\circ(L)\setminus\{v\})
$$
induced on chains by $\xi\mapsto\cone^\circ(\xi)\setminus\{v\}$.
As a consequence we have a quasi-isomorphism
 $$ I_\p C^\bullet(L)[n-k+1]\To
\mathcal I_\p\mathcal C^\bullet(V\cap U_k) =({\open_k}_*\mathcal
I_\p\mathcal C^\bullet)(V)
 $$
On the other hand Propositions~\ref{p:suspension} and \ref{p:coning}
give a quasi-isomorphism
 $$
\tau_{\le -n+\p(k)}(I_\p C^\bullet(L) [n-k+1])\To
\mathcal I_\p\mathcal C^\bullet(V).
$$
The resulting diagram
\begin{equation*}
\begin{array}{ccc}
\tau_{\le -n+\p(k)}(I_\p C^\bullet(L) [n-k+1])&\To
&\tau_{\le -n+\p(k)}({\open_k}_*\mathcal I_\p\mathcal C^\bullet)(V)\\
\downarrow& & \downarrow\\
\IpCsheaf^\bullet(V)&\overset{\alpha_k(V)}\To & \tau_{\le
-n+\p(k)}(R{\open_k}_*\mathcal I_\p\mathcal C^\bullet)(V)
\end{array}
\end{equation*}
commutes. Hence $\alpha_k$ is an isomorphism in $\mathcal
D^b(U_{k+1})$.
\end{proof}

\begin{cor}
The complexes $\mathcal I_\p\mathcal C^\bullet_X$ and $\mathcal
P_{\p}^\bullet$ are isomorphic in $\mathcal D^b(X)$. In particular we have
 $$
 I_\p H^*(X)=\mathbb H^*(\mathcal P_{\p}^\bullet(\mathbb O_{U_2})),
  $$
and intersection cohomology is independent of the choice of
\pl-structure on $X$. \qed
\end{cor}

The definition of $I_\p H^*(X)$ can easily be generalized to have
coefficients in a local system $\mathcal L$ on $U_2$. For any
simplex $\Delta$ occurring in a $\p$-allowable chain the
intersection conditions guarantee that $\Delta\cap U_2$ contains
the interior $\Delta^\circ$ of $\Delta$, so that it makes sense to
take the coefficient of $\Delta$ in $\Gamma(\Delta^\circ,\mathcal L)$.
With this definition the corollary generalizes to
\begin{equation*}
I_\p H^*(X,\mathcal L)=\mathbb H^*(\mathcal
P_{\p}^\bullet(\mathcal L\otimes\mathbb O_{U_2})).
\end{equation*}
\begin{rem}\label{r:0IC}
As another easy application let us give a proof of
Proposition~\ref{p:extremalperv}.(2), that $I_\0
H_i(X)=H^{n-i}(X)$ for a normal, oriented pseudomanifold $X$.
Note that since
 $$
\mathcal P_\p^\bullet(\mathcal L)= \left(\tau_{\le \p(n)}\circ R
{\open_n}_*\circ\cdots\circ \tau_{\le \p(3)} \circ R
{\open_3}_*\circ\tau_{\le\p(2)}\circ R {\open_2}_* (\mathcal
L)\right)[n]
 $$
and $\tau_{\le 0}\circ R(\open_m)_*=(\open_m)_*$ on $\Sh(U_m)$, we
have
 $$
 \mathcal P_\0^\bullet(\mathcal L)=\open_*(\mathcal L)[n],
 $$
where $\open :U_2\to X$ is the inclusion into $X$. So $\mathcal
P_\0^\bullet(\mathbb O_{U_2})=\open_*(\mathbb O_{U_2})[n]$. Since
$U_2$ is orientable, $\mathbb O_{U_2}=\Const_{U_2}$. Normality of
$X$ implies that $\open_*(\Const_{U_2})=\Const_X$ (since
intersecting open sets with $U_2$ preserves the number of
connected components). Therefore
 $$
I_\0 H_i(X)=\mathbb H^{-i}(\mathcal P^\bullet_\0)=\mathbb H^{-i}(\Const_X[n])=H^{n-i}(X).
 $$
\end{rem}

\section{Lecture - Verdier duality}

\subsection{The dual of a local system}
Recall that a local system $\mathcal L$ on $X$ (connected)
corresponds to a representation $V$ of the fundamental group, by
monodromy. The dual local system $\mathcal L^\vee$ has a natural
definition as the local system with monodromy representation
$V^*$ dual to $V$. Equivalently $\mathcal L^\vee$ may be given by
\begin{equation}\label{e:dualsheaf}
\mathcal L^\vee(U)=\Hom_{\Sh(U)}(\mathcal L|_U,\Const_U).
\end{equation}
For example $\Const_X^\vee=\Const_X$.

This simple definition of duality is just the right one for
orientable manifolds $M$ in that we have the following version
of  Poincar\'e duality:
 $$
H_c^i(M,\mathcal L)^*\cong H^{m-i}(M,\mathcal L^\vee).
 $$
This isomorphism comes from the usual Poincar\'e duality pairing
just with coefficients in taken in the dual local systems
$\mathcal L$ and $\mathcal L^\vee$, respectively.

\subsection{Duality for singular spaces}
Let $X$ be a singular pseudomanifold. The origin of intersection
cohomology's Poincar\'e duality on the level of $\mathcal D^b(X)$
is the duality functor of Borel and Moore \cite{BM:homology} and
J.-L.~Verdier \cite{Verdier:SemBour}. In rough outline, following
Verdier's approach, this duality takes the following form. Given
$X$ one constructs a {\it dualizing complex} in $\mathcal D^b(X)$,
which will turn out to be an injective resolution of the complex
of Borel-Moore chains $\mathcal D^\bullet_X$. Let us call it
$\dcomplex_X$. Then the dual of a complex $\mathcal
F^\bullet\in\mathcal D^b(X)$ is given by
 $$ \mathbb D_X(\mathcal
F^\bullet)=\mathcal Hom^\bullet(\mathcal F^\bullet,\dcomplex_X),
 $$
with notation $\mathcal Hom^\bullet$ defined in
Definition~\ref{d:Homdot} below.

\subsection{Versions of the Hom-functor}
Let us first define a sheaf from $\Hom$ between sheaves,
following the example of the dual of a local system
\eqref{e:dualsheaf}.
\begin{defn}\label{d:InternalHom}
For $\mathcal F,\mathcal G\in \Sh(X)$, define a $\mathcal
Hom(\mathcal F,\mathcal G)\in\Sh(X)$ by
 $$
\mathcal Hom(\mathcal F,\mathcal G)(U)=\Hom_{\Sh(U)}(\mathcal F|_U,\mathcal G|_U)
 $$
with the obvious restriction maps. Homomorphisms between sheaves
are locally determined, so this is indeed a sheaf. Note that
\begin{equation}\label{e:gammahom}
\Gamma\circ \mathcal Hom=\Hom_{\Sh(X)}.
\end{equation}
\end{defn}

\begin{defn}\label{d:Homdot}  For two bounded complexes
$\mathcal A^\bullet, \mathcal B^\bullet$ of sheaves let $\mathcal
Hom^\bullet(\mathcal A^\bullet,\mathcal B^\bullet)$ be the
bounded complex associated to the double complex $\mathcal
Hom(\mathcal A^i,\mathcal B^j )$. Explicitly,
 $$
\mathcal Hom^k(\mathcal A^\bullet,\mathcal
B^\bullet):=\prod_{i\in\Z}\mathcal Hom(\mathcal A^i,\mathcal
B^{i+k}),
 $$
where $d^k:\mathcal Hom^k(\mathcal A^\bullet,\mathcal
B^\bullet)\to \mathcal Hom^{k+1}(\mathcal A^\bullet,\mathcal
B^\bullet)$ has $\mathcal Hom(\mathcal A^i,\mathcal
B^{i+k+1})$-component given by $d^k(\phi)_i= d_{\mathcal
B}^{i+k}\circ \phi_k + (-1)^{i+1}\phi_{k+1} \circ d_{\mathcal
A}^i $.

In order to get a well defined bifunctor on the derived category,
say $\mathcal A^\bullet,\mathcal B^\bullet\in\mathcal D^b(X)$, we
need to replace $\mathcal B^\bullet$ by its (bounded) injective
resolution $I(\mathcal B^\bullet)$. So
\begin{equation*}
R\mathcal Hom^\bullet : \mathcal D^b(X)\x\mathcal
D^b(X)\to\mathcal D^b(X),
\end{equation*}
is defined by $R\mathcal Hom^\bullet(\mathcal A^\bullet,\mathcal
B^\bullet)= \mathcal Hom^\bullet(\mathcal A^\bullet,I(\mathcal
B^\bullet))$.
\end{defn}

\begin{defn} The global version
\begin{equation*}
\RHom^\bullet : \mathcal
D^b(X)\x\mathcal D^b(X)\to\mathcal D^b\Vect,
\end{equation*}
may be defined by $\Rhom^\bullet=R\Gamma\circ R\mathcal
Hom^\bullet$.
\end{defn}
The usual $\Hom$ in the derived category is recovered by
\begin{equation*}
\mathbb H^i(R\mathcal Hom^\bullet(\mathcal F^\bullet,\mathcal
G^\bullet))= \Hom_{\mathcal D^b(X)}(\mathcal F^\bullet,\mathcal
G^\bullet[i]).
\end{equation*}

\begin{defn}[Tensor products] For two sheaves $\mathcal F,\mathcal
G$ on $X$ their tensor product $\mathcal F\otimes\mathcal G$ is
defined as the sheafification of the presheaf $U\mapsto \mathcal
F(U)\otimes\mathcal G(U)$. Since we are dealing with sheaves of
vector spaces there are no difficulties with flatness and
tensoring with a given sheaf is an exact functor.
Setting $ (\mathcal A^\bullet\otimes\mathcal
B^\bullet)^k:=\bigoplus_{i+j=k}\mathcal A^i\otimes \mathcal B^j $
with the appropriate differential gives a well-defined bifunctor
$\otimes:\mathcal D^b(X)\x\mathcal D^b(X)\to\mathcal D^b(X)$.
\end{defn}

\begin{rem}
$\mathcal Hom(\mathcal
F,\underline{\ \, }\,)$ is right adjoint to the exact functor $
\mathcal F\otimes\underline{\ \,}\,$.
This implies that  $\mathcal Hom(\mathcal F,\underline{\ \, }\,)$ takes
injectives to injectives.
\end{rem}

\subsection{Pull-back with compact support}
We will follow the longer scenic route to constructing the
dualizing complex via the introduction of a new pull-back functor.
Let $f: X\to Y$ be a continuous map. Recall the two notions of
push-forward, $Rf_*$ and $Rf_!$. The usual pull-back functor $f^*$
is related to $Rf_*$ by the adjunction,
\begin{equation*}\Hom_{\mathcal D^b(X)}(f^*\mathcal
Y^\bullet,\mathcal X^\bullet)\cong
\Hom_{\mathcal D^b(Y)}(\mathcal Y^\bullet,Rf_*\mathcal
X^\bullet).
\end{equation*}

Duality should interchange cohomology with cohomology with proper
supports and thus, more generally, $Rf_*$ with $Rf_!$. Assuming
that, there should also be a dual notion of pull-back and a dual
adjunction formula. So we would like to construct a functor
$f^!:\mathcal D^b(Y)\to\mathcal D^b(X)$ such that
\begin{equation}\label{e:adjunction1} \Hom_{\mathcal D^b(X)}(\mathcal
X^\bullet,f^!\mathcal Y^\bullet)\cong
\Hom_{\mathcal D^b(Y)}(Rf_!\mathcal X^\bullet,\mathcal
Y^\bullet).
\end{equation}

\subsubsection{The case of an open embedding}\label{s:OpenInclusions}
If $\open :U\hookrightarrow X$ is an open inclusion then $\open^*$
is already right adjoint to $\open_!$,
 $$
\Hom_{\Sh(X)}(\mathcal X,\open^*\mathcal Y)\cong
\Hom_{\Sh(Y)}(\open_!\mathcal X,\mathcal Y),
 $$
with obvious adjunction morphisms $\mathcal X\to \open^*
\open_!\mathcal X$ and $\open_! \open^*\mathcal Y\to\mathcal Y$.
In this case we can simply set $\open^!\mathcal
F^\bullet=\open^*\mathcal F^\bullet$. Note that since $\open^*$
and $\open_!$ are both exact, the adjunction formula implies that
$\open^*$ takes injectives to injectives and $\open_!$ takes
projectives to projectives.

\subsubsection{The case of a closed embedding}
Let us first introduce some simplifying notation. Suppose
$\mathcal F\in\Sh(X)$ and $\open_U:U\hookrightarrow X$ is an open
inclusion. Define
 $$
\mathcal F_{U,X}:=(\open_U)_!(\mathcal F|_U) \in \Sh(X).
 $$
Explicitly, $\mathcal F_{U,X}(W)=\{s\in\mathcal F(U\cap W)\ |\
\text{$\Supp(s)$ is closed in $W$}\}$. For example
$\Const_{U,X}(W)=\C^c$ where $c$ is the number of connected
components of $W$ wholly contained in $U$.

Suppose $f=\closed:Z\to Y$ is the inclusion of a closed subset.
Define a functor $i':\Sh(Y)\to\Sh(Z)$ by
\begin{equation*}
i'(\mathcal Y)(U):=\Hom_{\Sh(Y)}(\closed_!\Const_{U,Z},\mathcal
Y)= \{\text{Sections $s\in\mathcal Y(V)$ supported in $Z$}\},
\end{equation*}
where $V$ is an open set in $Y$ such that $U=V\cap Z$. Then
$i'(\mathcal Y)$ really is a sheaf, since the property of being
supported in $Z$ can be checked locally. And $i'$ is clearly left
exact. Notice that this definition was engineered to give an
adjoint to $\closed_!$,
\begin{equation*}
\Hom_{\Sh(Z)}(\Const_{U,Z},i'(\mathcal Y))= i'(\mathcal
Y)(U)=\Hom_{\Sh(Y)}(\closed_!\Const_{U,Z},\mathcal Y).
\end{equation*}
The pull-back with compact support can be defined by
$\closed^!=Ri'$.

\subsubsection{Construction of $f^!$ for general $f$} Let
$f:X\to Y$ be an arbitrary continuous map. Then the presheaf
defined by
\begin{equation}\label{e:wrongdef}
U\mapsto\Hom_{\Sh(Y)}(f_! \Const_{U,X},\mathcal Y)
\end{equation}
is not necessarily a sheaf. For example it may have nontrivial
global sections while all of its stalks are zero. This is to do
with the inflexibility of the sheaf $\Const_{U,X}$ along the
fibers of $f$, which rarely allows for sections with proper
support. This problem is solved by going over to the derived
category and replacing $\Const_X$ with a fixed soft resolution
$\mathcal S_X^\bullet$. For example we may take the soft
resolution $\mathcal C_X^\bullet$ of $\Const_X$ from
Section~\ref{s:ICsheaves} if $X$ is piecewise linear.
\begin{thm}[{\cite[VI 1.1]{Iversen:SheafCohBook}}, {\cite[VI]{Borel:ICBook}}]
For $f:X\to Y$, the assignment
 $$
f^!(\mathcal Y^\bullet)(U):= \RHom^\bullet(Rf_!\,\mathcal
S_{U,X}^\bullet,\mathcal Y^\bullet)
 $$
defines a functor $f^!:\mathcal D^b(Y)\to\mathcal D^b(X)$. The
functor $f^!$ is (up to isomorphism) independent of the choice of
soft resolution $\mathcal S_X^\bullet$.
\end{thm}
Note that the discussion around \eqref{e:wrongdef} implies that
the functor $f_!:\Sh(X)\to\Sh(Y)$ simply has no right adjoint in
general. But $f^!$ does give a right adjoint functor to $Rf_!$
directly in the derived category. This, in a local version, is
the main result of Verdier \cite{Verdier:SemBour}, see
\cite[VII.5 Theorem~5.2]{Iversen:SheafCohBook}.

\begin{thm}\label{t:LocalVerdierDuality}
Let $\mathcal X^\bullet\in\mathcal D^b(X)$ and $\mathcal Y^\bullet
\in\mathcal D^b(X)$. Then there is a natural isomorphism
\begin{equation}\label{e:Verdier}
R f_*\, R\mathcal Hom^\bullet(\mathcal X^\bullet,f^!\mathcal
Y^\bullet) \cong R\mathcal Hom^\bullet(Rf_!\mathcal
X^\bullet,\mathcal Y^\bullet).
\end{equation}
\end{thm}

The adjunction formula \eqref{e:adjunction1} is an immediate
consequence.
\begin{cor} Applying $\mathbb H^0$ to both sides of \eqref{e:Verdier} we obtain a
natural isomorphism
\begin{equation*}\label{e:adjunction2}
\Hom_{\mathcal D^b(X)}(\mathcal X^\bullet,f^!\mathcal
Y^\bullet)\cong \Hom_{\mathcal D^b(Y)}(Rf_!\mathcal
X^\bullet,\mathcal Y^\bullet).
\end{equation*}
\end{cor}

\subsection{The dualizing complex}\label{s:dcomplex}
For the bounded derived category of
sheaves on a point the dualizing object is just
$\Const_{\{ pt\}}$ viewed as a  complex concentrated in degree
zero.
\begin{defn} The {\it dualizing complex} $\dcomplex_X\in\mathcal D^b(X)$
is defined by
 $$
\dcomplex_X:=\pi^!(\Const_{\{ pt\}}),
 $$
where $\pi=\pi_X:X\to\{pt\}$ is the projection to a point.
\end{defn}
Let us compute the dualizing complex from the definition, taking
as the soft resolution $\mathcal S_X^\bullet$ of $\Const_X$ the
complex of geometric cochains $\mathcal C_X^\bullet$. We get
\begin{align}\label{e:DoubleDual}
\pi^!\Const_{\{\pt\}}(U)=&\Hom^\bullet(\pi_!\mathcal C_{U,X}^\bullet,\C[0])=
\Hom^\bullet(\Gamma_c(\mathcal C_{U,X}^\bullet),\C[0])=
\Gamma_c(U,\mathcal C_X^{-\bullet})^*,
\end{align}
where we used that $\Const_{\{\pt\}}=\C$ is already injective.
This recovers the dual of the constant sheaf in the sense of
Borel and Moore~\cite{BM:homology}.

\subsection{Comparison with $\mathcal D^\bullet_X$ and other properties}
The complex $\mathcal D^\bullet_X$ of Borel-Moore chains is
related to the dualizing complex $\dcomplex_X$ by the following
map. From the definition of $\mathcal C_X^{i}$,
Section~\ref{s:ICsheaves}, one can see that there is a dual
pairing $\Gamma_c(U,\mathcal C_X^{-i})\x \Gamma(U,\mathcal
D^{i}_X)\to \C$. Therefore we have an inclusion
 $$
\mathcal D_X^{i}(U)\hookrightarrow\Gamma_c(U,\mathcal C^{-i}_X)^*
 $$
of $\mathcal D_X^{i}(U)$ into its double dual. Since these
complexes have finite dimensional cohomology, the inclusion
induces isomorphisms there. So we get a quasi-isomorphism
$\mathcal D^\bullet_X\to\dcomplex_X$.

One major advantage the dualizing complex has over the complex of
Borel-Moore chains is that it is a complex of injectives, while
$\mathcal D^\bullet_X$ was only soft. To see this intuitively,
imagine for the sections of $\Gamma_c(U,\mathcal
C_X^{-\bullet})^*$ something like geometric chains but no longer
necessarily with closed support. Then every restriction map is
split by an extension by zero, and the corresponding sheaf is
injective by Lemma~\ref{l:Spaltenstein}. A proper, more general
proof shows that there is a natural isomorphism
 $
 \Hom_{\Sh(X)}(\mathcal A ,\hat\dsheaf_X^{-i})=\Gamma_c(X,\mathcal
A\otimes\mathcal S^i_X)^*
 $
generalizing \eqref{e:DoubleDual} (e.g.
\cite[V Proposition 1.5]{Iversen:SheafCohBook}). Then since
$\mathcal S^i_X$ is soft, and flat of course,  $\mathcal
A\mapsto\mathcal A\otimes\mathcal S^i_X$ is an exact functor
taking $\mathcal A$ to a soft sheaf. Since soft sheaves are
$\Gamma_c$-acyclic, the combined functor
$\Hom_{\Sh(X)}(\,\underline {\ } \, ,\hat\dsheaf_X^{-i})$ is
exact.

Another advantage of $\dcomplex_X$ over $\mathcal D^\bullet_X$ is
that it can be defined in more general settings. We may choose any
soft resolution of $\Const_X$ to define $\pi^!$. In particular,
$X$ need not have a piecewise linear structure.

Note also that the dualizing complexes are inherently compatible
with pull-back with compact supports.  For a composition $f\circ
g$ of two continuous maps there is a natural isomorphism $(f\circ
g)^!= g^!\circ f^!$. And therefore
 $$
f^!(\dcomplex_Y)= f^!\circ\pi^!_Y(\C)=(\pi_Y\circ
f)^!(\C)=\pi_X^!(\C)=\hat{\mathcal D}^\bullet_X.
 $$
This is analogous to the property $f^*(\Const_Y)=\Const_X$.

\subsection{Verdier duality}\label{s:VDProperties}
The most important property of $\dcomplex_X$ is
that it is a dualizing complex in the categorical sense
for the duality for (constructible) complexes of
sheaves introduced by Borel and Moore.

\begin{defn} For $\mathcal X^\bullet\in\mathcal D^b(X)$ let
\begin{equation*}
\mathbb D_X(\mathcal X^\bullet):=\mathcal RHom^\bullet(\mathcal
X^\bullet, \dcomplex_X).
\end{equation*}
This defines a contravariant functor
$\mathbb D_X:\mathcal D^b(X)\to\mathcal D^b(X)$
which is called the {\it Verdier duality functor}.
\end{defn}

Verdier's main theorem, the local adjunction formula \eqref{e:Verdier},
now implies
\begin{equation}\label{e:lVD}
\mathbb D_X(\mathcal X^\bullet)=\mathcal
RHom^\bullet(R\pi_!(\mathcal X^\bullet),\Const_{\{\pt\}}),
\end{equation}
where $f=\pi:X\to\{\pt\}$.
The right hand side of \eqref{e:lVD} agrees with the
definition of a dual of Borel and Moore \cite{BM:homology}.
Let us rewrite $R\pi_!$ as $R\Gamma_c$ and also apply
$R\pi_*=R\Gamma$ to \eqref{e:lVD}. So
\begin{equation*}
R\Gamma(\mathbb D_X(\mathcal X^\bullet))=
R\Hom^\bullet(R\Gamma_c(\mathcal X^\bullet),\Const_{\{\pt\}})
\end{equation*}
in $\mathcal D^b\Vect$. Then taking $H^{-i}$ on both sides gives
the very appealing global sense,
\begin{equation}\label{e:globalVerdier}
\mathbb H^{-i}(X,\mathbb D_X(\mathcal X^\bullet))=\mathbb
H^{i}_c(X,\mathcal X^\bullet)^*,
\end{equation}
in which  $\mathbb D_X(\mathcal X^\bullet)$ is  dual to $\mathcal
X^\bullet$.

If for example $\mathcal X^\bullet=\Const_X$, then $\mathbb
D_X(\Const_X[0])= \dcomplex_X\cong\dsheaf_X^\bullet$. Therefore
\eqref{e:globalVerdier} becomes
\begin{equation*}
\begin{array}{ccc}
\mathbb H^{-i}(X,\mathcal D_X^\bullet)&=&\mathbb H^i_c(X,\Const_X)^*\\
 || &&|| \\
 H_i(X)&=&H_c^i(X)^*,
\end{array}
\end{equation*}
which is a special case of the universal coefficient theorem for cohomology.

\begin{rem}\label{r:VDProperties}
 Here are some basic properties of Verdier duality.
\begin{enumerate}
\item $\mathbb D_X:\mathcal D^b(X)\to\mathcal D^b(X)$ is a contravariant
functor. With respect to shifts, $\mathbb D_X(\mathcal
F^\bullet[n])=\mathbb D_X (\mathcal F^\bullet)[-n]$.
\item Let $f:X\to Y$. Then
\begin{align}
Rf_*\circ\mathbb D_X&=\mathbb D_Y\circ Rf_!,\\
\mathbb D_X\circ f^* &=f^!\circ \mathbb D_Y.\label{e:fupper!D}
\end{align}
In particular, Verdier duality intertwines $R\Gamma=R\pi_*$ and
$R\Gamma_c=R\pi_!$. Here the first formula follows from the
adjunction \eqref{e:Verdier}. For the second one see
\cite[V \S10]{Borel:ICBook}.
\item Suppose $\mathcal L$ is a local system on an $m$-dimensional
manifold $M$.
Then $\dcomplex_M\cong\mathbb O_M[m]$ and, since $\mathcal L$ is projective,
$R\mathcal Hom^\bullet(\mathcal L[m],\dcomplex_X)=\mathcal
Hom^\bullet(\mathcal L[m],\mathbb O_M[m])$. Therefore
$$
\mathbb D_M(\mathcal L[m])=\mathcal L^\vee\otimes\mathbb O_M.
$$
\end{enumerate}
\end{rem}

\subsection{Constructibility} For a local system on a manifold
$M$ the final property in Remark~\ref{r:VDProperties} implies that
$\mathbb D_M(\mathbb D_M(\mathcal L))\cong\mathcal L$. This has
an analogue for stratified spaces $X$.

\begin{defn}\label{d:constructible}
Let $\sigma$ be a pseudomanifold stratification of $X$. A complex
$\mathcal F^\bullet\in\mathcal D^b(X)$ is called {\it
$\sigma$-constructible} if for each stratum $S$ we have that
$\mathcal H^i(\mathcal F^\bullet|_{S})$ is a local system.
$\mathcal F^\bullet$ is called {\it constructible} if it is
$\sigma$-constructible for some pseudomanifold stratification
$\sigma$ of $X$.
\end{defn}

Note that all of the complexes we have been studying are
constructible for the (given) stratification on $X$. For example
the Deligne complex $\mathcal P^\bullet$ is constructible. This is
essentially because we `constructed' it from the constant sheaf,
or a local system, using only truncation functors and
push-forwards $R(\open_k)_*$, where the inclusions $\open_k$ were
compatible with the stratification (and the stratification obeys
local normal triviality).

The dualizing complex $\dcomplex_X$ is in fact constructible for
any pseudomanifold stratification $\sigma$ of $X$. For example if
$U\subset X$ open is a manifold, then we know that $$
\dcomplex_X|_U\cong\mathbb O_U[n] $$ by calculation of the local
Borel-Moore homology. Similarly one can deduce that $\dcomplex_X$
has finite-dimensional locally constant cohomology along the
strata of $X$, using local normal triviality. Moreover, if
$\mathcal F^\bullet$ is a $\si$-constructible complex, then so is
$\mathbb D_X(\mathcal F^\bullet)$. See
\cite[V \S8]{Borel:ICBook}.

Finally, if $\mathcal F^\bullet\in\mathcal D^b(X)$ is
constructible, then the natural map
 $$\mathcal F^\bullet\to\mathbb
D_X\circ\mathbb D_X(\mathcal F^\bullet)
 $$
is an isomorphism in $\mathcal D^b(X)$. See 
\cite[V Theorem 8.10]{Borel:ICBook}.

\section{Lecture - Topological invariance of intersection cohomology and
Poincar\'e duality}

The main aim of this section will be to apply the Verdier duality
functor to the Deligne sheaf.

\subsection{Long exact sequences}
Suppose $\mathcal F\in\Sh(X)$ and $X=Z\sqcup U$ is a
decomposition of $X$ into a closed subset $Z$ and an open one
$U$. Denote the inclusions $\closed:Z\hookrightarrow X$ and $\open
:U\hookrightarrow X$. Let $\Gamma_Z(\mathcal F):=\closed_!
i'(\mathcal F)\in\Sh(X)$. So
 $$
\Gamma_Z(\mathcal F)(V)=\Gamma_Z(V,\mathcal F)=\{s\in\mathcal
F(V)\ |\ \Supp(s)\subset Z\}.
 $$
Then the decomposition of $X$ into $Z$ and $U$ gives rise to an
exact sequence
\begin{equation*}
0\to\Gamma_Z(\mathcal F)\to\mathcal F\to \open_* \open^*\mathcal F
\end{equation*}
which extends to a short exact sequence $\to 0$ if $\mathcal F$ is injective.

Suppose now $\mathcal F^\bullet$ is a complex, and $\mathcal
I^\bullet=I(\mathcal F^\bullet)\cong \mathcal F^\bullet$ its
injective resolution. Then we get a sequence of chain maps
\begin{equation*}
0\to \closed_! \closed^!\mathcal I^\bullet\to\mathcal
I^\bullet\to \open_* \open^*\mathcal I^\bullet\to 0
\end{equation*}
which is exact in every degree. In the language of triangulated
categories, we obtain a {`distinguished triangle'},
\begin{equation}\label{e:triangle1}
\begin{array}{ccc}
\closed_!\closed^!\mathcal F^\bullet&\To &\quad\mathcal F^\bullet\\
\quad [1]\nwarrow & &  \swarrow\quad \\
&R\open_*\open^*\mathcal F^\bullet
\end{array},
\end{equation}
where the map on the left hand side is a map to the shifted
complex $R\closed_!\closed^!(\mathcal F^\bullet)[1]$. Applying
cohomology functors such as $\mathbb H^\bullet$, $\mathbb
H^\bullet_c$, or $\mathcal H^\bullet_x$ gives a long exact
sequence. For example,
\begin{equation*}
\cdots\to\mathbb H^i(Z,\closed^!\mathcal F^\bullet)\to\mathbb
H^i(X,\mathcal F^\bullet)\to \mathbb H^i(U,\mathcal F^\bullet)\to
\mathbb H^{i+1}(Z,\closed^!\mathcal F^\bullet)\to \cdots.
\end{equation*}

The inclusions $i$ and $j$ also define another distinguished
triangle,
 $$
\cdots\to j_! j^!\mathcal F^\bullet \to \mathcal F^\bullet \to i_* i^*\mathcal F^\bullet \overset{[1]}\to\cdots .
 $$

\subsection{Costalks}
Let $\closed_x:\{x\}\to X$ be the inclusion of a point. The dual
concept to the stalk functor $\closed_x^*$ is the functor of
costalks
$\closed_x^!:\mathcal D^b(X)\to\mathcal D^b(\{pt\})$.
Define $$ \mathcal H^j_{\{x\}}(\mathcal
F^\bullet):=H^j(\closed_x^!\mathcal F^\bullet).$$ This is called
the local cohomology supported in $\{x\}$.

Note that the stalk cohomology of $\mathcal F^\bullet$ at $x$ can
be described as
\begin{equation*}
\mathcal H^j_{x}(\mathcal F^\bullet)=\lim_{\overset\To{\{U \text{
open}\, |\, x\in U\}}} \mathbb H^j(U,\mathcal F^\bullet),
\end{equation*}
by replacing $\mathcal F^\bullet$ by its injective resolution in
\eqref{e:stalks}. If $\mathcal F^\bullet$ is constructible,
then an analogous description,
\begin{equation*}
\mathcal H^j_{\{x\}}(\mathcal
F^\bullet)=\lim_{\overset{\longleftarrow}{\{U \text{ open}\, |\,
x\in U\}}} \mathbb H^j_c(U,\mathcal F^\bullet),
\end{equation*}
for the costalk follows by Verdier duality.  See also
\cite[V \S3]{Borel:ICBook}.

\begin{lem}\label{l:costalks}
Suppose $M$ is an $m$-manifold and $\mathcal F^\bullet\in\mathcal
D^b(M)$ has finite rank, locally constant cohomology. Let $x\in M$
and $U\subset M$ an open ball containing $x$. Then
\begin{itemize}
\item[$(a)$] $\mathcal H_x^j(\mathcal F^\bullet)=\mathbb H^j(U,\mathcal F^\bullet)$.
\item[$(b)$] $\mathcal H^j_{\{x\}}(\mathcal F^\bullet) =\mathbb H^j_c(U,\mathcal F^\bullet)
=\mathcal H^{j-m}_x(\mathcal F^\bullet).$
\end{itemize}
In particular one has $j_x^!(\mathcal F^\bullet)\cong
j_x^*(\mathcal F^\bullet)[-m]$.
\end{lem}
\begin{proof}
To see $(a)$, consider the spectral sequence
\begin{equation*}
E^{ij}_2:H^i(U,\mathcal H^j\mathcal F^\bullet|_U)
\Rightarrow \mathbb H^{i+j}(U,\mathcal F^\bullet).
\end{equation*}
Since $\mathcal H^j(\mathcal F^\bullet|_U)$ is a constant sheaf,
by assumption, the cohomology $H^i(U,\mathcal H^j\mathcal
F^\bullet|_U)$ is trivial for $i>0$. Therefore the spectral
sequence collapses to give $\mathbb H^j(U,\mathcal F^\bullet)=
H^0(U,\mathcal H^j(\mathcal F^\bullet))=\mathcal H^j_x(\mathcal
F^\bullet)$. The second equality in $(b)$ is proved in a similar
way.

To get the identity $\mathcal H^j_{\{x\}}(\mathcal F^\bullet)
=\mathbb H^j_c(U,\mathcal F^\bullet)$ apply Verdier duality to
$(a)$, in particular \eqref{e:globalVerdier} and
\eqref{e:fupper!D}.
\end{proof}

\subsection{Costalks and the Deligne sheaf}\label{s:DeligneCostalks}
As usual let $X$ be an $n$-dimensional stratified pseudomanifold
with a local system $\mathcal L$ on the open stratum
$U_2=X\setminus\Sigma$, and $\p$ a fixed perversity. We use the
same notation as in Section~\ref{s:towardsDeligne}. The Deligne
sheaf $\mathcal P^\bullet= \mathcal P^\bullet_{\p}(\mathcal L)$
in $\mathcal D^b(X)$ has  the following properties (which also
characterize it up to isomorphism).
\begin{enumerate}
\item[$(\, A1\, )$] $\mathcal P^\bullet$ is constructible for the given stratification on $X$.
\item[$(\, A2\, )$]
$\mathcal P^\bullet|_{U_2}\cong\mathcal L[n]$.
\item[$(\, A3\, )$] $\mathcal H^i_x(\mathcal P^\bullet)=0$ for $i>-n+\p(k)$ and $x\in S_k$.
\item[$(\, A4\, )$] Consider the inclusion
$U_k\overset{\open_k}\hookrightarrow U_{k+1}$. The attachment map
 $$\alpha_k:\mathcal P^\bullet|_{U_{k+1}}\to R(\open_k)_*(\mathcal P^\bullet|_{U_k}) $$
is a quasi-isomorphism up to degree $\le -n+\p(k)$.
\end{enumerate}

\vskip .2cm

Let us prove that $\mathcal P^\bullet$ also satisfies the
following `dual' property to $(A3)$~:
\begin{itemize}
\item[$(A3^!)$] $\mathcal H^i_{\{x\}}(\mathcal P^\bullet)=0$ for $x\in S_k$ and  $i<-\underline
q(k)$,
where $\underline q$ is the complementary perversity to $\p$.
\end{itemize}
\begin{lem}\label{l:P5}
$\mathcal P^\bullet$ satisfies $(A3^!)$
\end{lem}

\begin{proof}
Let $\mathcal P^\bullet_k:=\mathcal P^\bullet|_{U_k}$ and
consider the inclusions $$ U_k\overset{j}\hookrightarrow
U_{k+1}\overset i\hookleftarrow S_k. $$ Let $x\in S_k$. Applying
$\mathcal H^l_x$ to the distinguished triangle \eqref{e:triangle1} gives a
long exact sequence
\begin{equation*}
\cdots\to \mathcal H^l_x(\closed^!\mathcal P^\bullet_{k+1})\to
\mathcal H^l_x(\mathcal P^\bullet)\to \mathcal H^l_x(R \open_*
\open^*\mathcal P^\bullet_{k+1}) \to\mathcal
H^{l+1}_x(\closed^!\mathcal P^\bullet_{k+1})\to\cdots
\end{equation*}
To use this long exact sequence, we need to relate the terms
$\mathcal H^l_x(\closed^!\mathcal P^\bullet_{k+1})$ to costalks.
Let $\closed_x:\{x\}\to S_k$ denote the inclusion. Then one can
use Lemma~\ref{l:costalks} to deduce
\begin{align*}
\mathcal H^l_x(\closed^!\mathcal
P^\bullet_{k+1})&=H^{l}(\closed_x^*\closed^!\mathcal
P_{k+1}^\bullet)=
H^{l+n-k}(\closed_x^!\closed^!\mathcal P_{k+1}^\bullet)\\
&=H^{l+n-k}((i\circ \closed_x)^!\mathcal P_{k+1}^\bullet)=
\mathcal H^{l+n-k}_{\{x\}}(\mathcal P_{k+1}^\bullet).
\end{align*}

Recall that by property $(A3)$ we have $\mathcal H^{l}_x(\mathcal
P^\bullet)=0$ for all $l>m:=-n+\p(k)$. So let us write down again
the long exact sequence,
\begin{equation}\label{e:costalkLES}
\begin{array}{ccccccc}
&&\cdots &\to &\mathcal H_x^{m-1}(\mathcal P^\bullet)&\to &
 \mathcal H^{m-1}_x(R \open_* \open^*\mathcal P^\bullet_{k+1})\\
& \to &
\mathcal H^{m+n-k}_{\{x\}}(\mathcal P^\bullet_{k+1})&\to &
\mathcal H^m_x(\mathcal P^\bullet)&\to &
\mathcal H^m_x(R \open_* \open^*\mathcal P^\bullet_{k+1})\\
&\to &\mathcal H^{m+n-k+1}_{\{x\}}(\mathcal P^\bullet_{k+1})&
\to &\mathcal H^{m+1}_x(\mathcal P^\bullet)&=0.
\end{array}
\end{equation}
From the property $(A4)$ of the attachment map we have that the
maps $\mathcal H_x^{l}(\mathcal P^\bullet)\to \mathcal H^{l}_x(R
\open_* \open^*\mathcal P^\bullet_{k+1})$ are isomorphisms
exactly for $l\le m$. But therefore the long exact sequence
\eqref{e:costalkLES} implies that
\begin{align*}
\mathcal H^{l+n-k}_{\{x\}}(\mathcal P^\bullet_{k+1})=0&\text{ for
$l< -n+\p(k)+2$}.
\end{align*}
Since $\underline q(k)=k-2-\p(k)$, this is precisely the condition
$(A3^!)$.
\end{proof}

\begin{prop}\label{p:A1-5} The Deligne sheaf
$\mathcal P^\bullet=\mathcal P^\bullet_\p(\mathcal L)\in \mathcal
D^b(X)$ is uniquely determined up to isomorphism by the properties
$(A1)$, $(A2)$, $(A3)$ and $(A3^!)$.
\end{prop}

This is proved basically by reading the proof of Lemma~\ref{l:P5}
in reverse. The vanishing of the costalks implies that the
attachment map is a quasi-isomorphism in the required degrees by
the same long exact sequence. So we have properties $(A1)-(A4)$.
But $(A3)$ and $(A4)$ together with the constructibility property
$(A1)$ imply that the attachment map gives a quasi-isomorphism
 $$
\mathcal P^\bullet_{k+1}\to\tau_{\le
-n+\p(k)}R(\open_{k})_*\mathcal P^\bullet_{k}.
 $$
So we recover the inductive construction of $\mathcal P^\bullet$.

\begin{thm}
Let $\mathcal L$ be a local system on the open stratum $U_2$ of 
an $n$-dimensional pseudomanifold $X$.
Let $\p$ and $\q$ be complementary perversities.
Then
 \begin{equation}\label{e:DofP}
\mathbb D_X(\mathcal P_\p^\bullet(\mathcal L))\cong
\mathcal P_\q^\bullet (\mathcal L^\vee\otimes \mathbb O_{U_2})[-n].
 \end{equation}
\end{thm}
This theorem is proved by using the properties of Verdier duality
to show that $\mathbb D_X(\mathcal P_\p^\bullet(\mathcal L))[n]$
satisfies $(A1),(A2),(A3)$ and $(A3^!)$ for $\mathcal
L^\vee\otimes\mathbb O_{U_2}$ and $\q$, and then applying
Proposition~\ref{p:A1-5}.

If we apply $\mathbb H^i$ to \eqref{e:DofP}, we obtain for the left hand side
\begin{equation*}
H^i\circ R\Gamma\circ\mathbb D_X(\mathcal P^{\bullet}_\p (\mathcal L))=
\mathbb H^{-i}_c(\mathcal P^\bullet_{\p}(\mathcal L))^*
\end{equation*}
using the property $R\Gamma\circ\mathbb D_X=\mathbb
D_{\{pt\}}\circ R\Gamma_c$ from Remark~\ref{r:VDProperties}. And
the right hand side becomes $\mathbb H^{-n+i}(\mathcal
P^\bullet_\q(\mathcal L^\vee\otimes\mathbb O_{U_2}))$.

Thus we obtain  Theorem~\ref{t:PD} as a corollary (in a slightly more general form).
\begin{cor}[Poincar\'e Duality]
With notation as above,
\begin{equation*}
I_\p H_c^{-i}(X,\mathcal L\otimes\mathbb O_{U_2} )^*\cong I_\q
H^{-n+i}(X,\mathcal L^\vee).
\end{equation*}
\end{cor}

\subsection{Variation of the stratification}
Let $X$ be a pseudomanifold without a fixed stratification. And
let $(U,\mathcal L)$ be a pair consisting of an open dense
submanifold $U\subset X$ with complement of codimension $\ge 2$
and a local system defined on $U$. We may assume $U$ is maximal
(by Zorn's Lemma applied to the set of extensions of $\mathcal
L$).

\begin{thm}[{\cite[V 4.15]{Borel:ICBook}}]\label{t:Ax2}
Suppose we are given a pair $(U,\mathcal L)$ as above, and assume that $U$ contains the open
stratum of some stratification of
$X$.\footnote{
Borel indicates an example of a pair
$(U,\mathcal L)$ where $U$ is maximal for $\mathcal L$ and does not contain the
open stratum $U_2$ of any stratification, 
see \cite[V Remark 4.14]{Borel:ICBook}.
However the $U$ in the example has no universal
cover (not locally simply connected), and the definition of local system used seems to be
as finite-dimensional representation of the fundamental group. For a local system in
the locally constant sheaf sense I don't know whether it may not
be possible to leave out this condition on $U$.
}
Let $\p$ be a perversity. Then there exists a unique (up to isomorphism) complex
$\tilde{\mathcal P}^\bullet\in\mathcal D^b(X)$ such that
\begin{itemize}
\item[$(\,B1\, )$] $\tilde{\mathcal P}^\bullet$ is constructible.
\item[$(\,B2\, )$] $\tilde{\mathcal P}^\bullet|_V\cong\mathcal L|_V[n]$ for some open dense
$V\subset U$ with codimension $\ge 2$ complement.
\item[$(\,B3\, )$]
$\dim(\Supp \mathcal H^{-n+i}(\tilde{\mathcal P}^\bullet))\le
n-\p\inv(i)$, where $\p\inv(i):=min\{k\, |\, \p(k)\ge i\}$.
\item[$(B3^!)$]
$\dim(\Cosupp \mathcal H^{-i}(\tilde{\mathcal P}^\bullet))\le
n-\q\inv(i)$, where the cosupport of a sheaf is defined by
$\Cosupp(\mathcal F):= \{x\in X |\mathcal F_{\{x\}}\ne 0\} $.
\end{itemize}
\end{thm}

Suppose $\si$ is a stratification of $X$ adapted to $\mathcal L$, that is,
with open stratum $U_2\subset U$. Then one can check that the
corresponding Deligne sheaf $\mathcal P^\bullet=\mathcal P_\p^\bullet(\mathcal
L)$ satisfies the above conditions. For example, property
$(A3)$ in Section~\ref{s:DeligneCostalks} immediately implies
that $\Supp( \mathcal H^{-n+i}({\mathcal P}^\bullet))$ cannot
contain any stratum $S_k$ having codimension less than $min\{k\,
|\, \p(k)\ge i\}$, so property $(B3)$ follows. For $(B3^!)$ it is
necessary to use also that Verdier duality preserves
$\sigma$-constructibility.

The idea of the proof of Theorem~\ref{t:Ax2} is to construct a
coarsest possible stratification $\sigma$ and
$\sigma$-constructible Deligne sheaf $\mathcal
P^\bullet_{\p,\sigma}$, simultaneously. For any other adapted
stratification $\sigma'$, the strata of $\sigma$ can be made up
of unions of strata of $\sigma'$. Then one can show that $\mathcal
P^\bullet_{\p,\sigma}$  also satisfies the axioms $(A1)-(A3^!)$
with regard to $\sigma'$ (using the $\sigma$-constructibility
property). Therefore by Proposition~\ref{p:A1-5} the two versions
of the Deligne sheaf must be isomorphic.

\begin{cor}
The Deligne sheaf $\mathcal P^\bullet_\p(\mathcal L)\in\mathcal
D^b(X) $ is independent of the stratification of $X$ up to
isomorphism. In particular, intersection cohomology is a
topological invariant.
\end{cor}

\subsection{Intersection cohomology of varieties and perverse sheaves}
Let $X$ be a quasi-projective algebraic variety over $\C$. Interesting examples
in the context of intersection cohomology include 
affine varieties such as the nilpotent
cone $\mathcal N$ in a Lie algebra, or 
projective varieties such as Grassmannians,
flag varieties and Schubert varieties.

Any quasi-projective variety $X$ has a stratification $X\supset
X_2\supset X_4\supset\cdots \supset X_{2n}$ by Zariski closed
subsets, such that the strata are smooth. A canonical such
stratification is obtained for example by taking the singular
locus of $X$ and the singular locus of the singular locus, and so
forth. A better stratification, a `Whitney stratification', which
can be shown also to satisfy local normal triviality can be
obtained by refinement (by a combination of \cite{Whitney:Strat}
and \cite[Part IV]{Borel:ICBook}, see also \cite[\S3.2]{Kirwan:IHBook}).
As a consequence $X$ has a pseudomanifold stratification by
closed subvarieties. See also \cite[Chapter 8]{KashScha:SheavesBook} and 
the nice introduction to Whitney stratifications in 
\cite[Part I Chapter 1]{GorMacPh:book}. In general 
references for the theory of perverse sheaves include 
\cite{BBD:PerverseSheaves}, \cite{KashScha:SheavesBook},
\cite{KiehlWeiss:WeilConjBook},
and the encyclopedic \cite{GelfandManin:HomologicalAlgebraBook}.

Let us restrict our attention to the middle perversity case, $\underline
m(2k)=k-1$.
Since $X$ and its strata are all even-dimensional there are
better grading conventions to which we now switch. Let
$n:=\dim_\C(X)=\frac 12\dim_\R(X)$. Then shift the
dualizing complex and Deligne sheaves
by $[-n]$. These should become
concentrated in degrees $-n\le
i\le n$, rather than between $-2n$ and $0$. The new version
$\mathcal P^\bullet= \mathcal P_{\underline m}^\bullet(\mathcal
L)\in\mathcal D^b(X)$ of the Deligne sheaf that can be
constructed in this setting has the following properties.
\begin{enumerate}
\item[$(\, C1\, )$]  $\mathcal P^\bullet$ is constructible for
a Whitney stratification of $X$ by Zariski closed subsets.
\item[$(\, C2\, )$] $\mathcal P^\bullet|_U\cong\mathcal L[\dim_\C(X)]$ for some Zariski open
dense $U\subset X$.
\item[$(\, C3\, )$] $\dim_\C (\Supp^{-i}\mathcal P^\bullet) <i$ for all $i<\dim_\C(X)$.
\item[$(C3^!)$] $\dim_\C (\Cosupp^{i}\mathcal P^\bullet) <i $ for all $i<\dim_\C(X)$.
\end{enumerate}

\begin{defn}
Let $\mathcal D^b_c(X)$ be the full subcategory of $\mathcal
D^b(X)$ of complexes with cohomology constructible for a Whitney
stratification of $X$ by Zarsiki closed subsets.

A complex $\mathcal G^\bullet\in\mathcal D^b_c(X)$ is called a
{\it perverse sheaf\,} if
\begin{align}\label{e:le0}
&\dim_\C(\Supp^{-i}\mathcal G^\bullet)\le i,\\
&\dim_\C(\Cosupp^{i}\mathcal G^\bullet)\le i, \label{e:ge0}
 \end{align}
for all $i\in\Z$. Let $\mathcal M(X)\subset\mathcal D^b_c(X)$ be
the full subcategory of perverse sheaves. If we fix a
stratification $\si$ of $X$ then we may also consider the
subcategory $\mathcal M_\si(X)$ of $\sigma$-constructible perverse
sheaves.
\end{defn}

This construction of a subcategory of $\mathcal D_c^b(X)$ is akin
to the description of $\Sh(X)$ inside $\mathcal D^b(X)$ as
$\mathcal D^{\le 0}(X)\cap \mathcal D^{\ge 0}(X)$, complexes of
sheaves concentrated in degree~$0$. Here \eqref{e:le0} and
\eqref{e:ge0} define the analogues of $\mathcal D^{\le 0}(X)$ and
$\mathcal D^{\ge 0}(X)$ --  they can be shown to satisfy the axioms of a
`$t$-structure' in the triangulated category. The resulting
category $\mathcal M(X)$ is the `core' of this $t$-structure and
therefore an abelian category by general theory.  We also have
that $\mathcal M(X)$ and $\mathcal M_\sigma(X)$ are preserved
under Verdier duality, thanks to the change of grading and the
choice of middle perversity. And the simple objects are known. In
summary~:

\begin{thm}[\cite{BBD:PerverseSheaves},
{\cite[Chapter III]{KiehlWeiss:WeilConjBook}}] Let $X$ be a
quasi-projective algebraic variety over $\C$ with a Whitney
stratification by Zariski closed subsets $\sigma$.
\begin{enumerate}
\item
The categories  $\mathcal M(X)$ and $\mathcal M_\si(X)$ are
abelian categories.
\item
The Verdier duality functor defines an anti-equivalence of
categories $\mathbb D_X:\mathcal D_c^b(X)\to \mathcal D^b_c(X)$
which restricts to give anti-equivalences $\mathcal M(X)\to
\mathcal M(X)$ and $\mathcal M_\si(X)\to \mathcal M_\si(X)$.
\item
For a pair $(S,\mathcal L)$ where $S\subset X$ is a closed
irreducible subvariety and $\mathcal L$ a local system on a smooth
Zariski open
subset of $S$ define
 $$\ICsheaf^\bullet_{
S,X}(\mathcal L):=\closed_!\mathcal P_{\underline
m}^\bullet(\mathcal L),
 $$
where $\closed:S\hookrightarrow X$ is the inclusion.  If
$\mathcal L$ is an irreducible local system, then
$\ICsheaf^\bullet_{ S,X}(\mathcal L)$ is a simple object in
$\mathcal M(X)$. All simple objects in $\mathcal M(X)$ are
isomorphic to $\ICsheaf^\bullet_{S,X}(\mathcal L)$ for some pair
$(S,\mathcal L)$ with $\mathcal L$ irreducible.
\end{enumerate}
\end{thm}

\begin{defn}
A complex $\mathcal F^\bullet\in\mathcal D^b(X)$ is called {\it semisimple}
if  it is of the form
 $$
\mathcal F^\bullet\cong\bigoplus_{(S,\mathcal L),k} \underline
V_{(S,\mathcal L)}^k\otimes\ICsheaf^\bullet_{S,X}(\mathcal L)[k]
 $$
for some finite dimensional $\Z$-graded multiplicity vector spaces
$V_{(S,\mathcal L)}=\oplus_{k\in\Z} V_{(S,\mathcal L)}^k$ which
are nonzero for only finitely many pairs $(S,\mathcal L)$.
\end{defn}

We may now state the famous theorem of Beilinson, Bernstein,
Deligne and Gabber \cite{BBD:PerverseSheaves}.
\begin{thm}[Decomposition theorem]
If $f:X\to Y$ is a proper morphism, then $Rf_*\Const_X$ is semisimple.
\end{thm}

Note that it is not known whether it is 
possible in the statement to replace the
constant sheaf by an arbitrary local system. The only allowable
local systems in this context are the ones that arise from the
constant sheaf on an algebraic covering. In other words, suppose
$p:\tilde X\to X$ is an unramified finite morphism. Then $p_*$ is
an exact functor, and $Rp_*\Const_{\tilde X}=p_*\Const_{\tilde
X}=\mathcal L$ is a local system on $X$ for which the
decomposition theorem (obviously) holds.

The original proof \cite{BBD:PerverseSheaves} of this theorem involves
base change to finite charactersitic and the 
study of eigenvalues of the Frobenius map on $\ell$-adic cohomology.
Another reference for the $\ell$-adic theory of perverse sheaves
(including the decomposition theorem) is
\cite{KiehlWeiss:WeilConjBook}. A proof of the decomposition theorem
without the base change was given by M. Saito
\cite{Saito:Decompthm}.  We now turn
to an application.

\section{Lecture - Kazhdan-Lusztig polynomials}

Let $W$ be a Weyl group, for example the symmetric group $S_n$.
Kazhdan and Lusztig \cite{KaLus:Hecke} introduced polynomials
$P_{u,w}(q)$ for pairs $u,w\in W$ and conjectured an
interpretation of their value at $1$ in terms of
decomposition numbers of Verma modules
in the representation theory of the
corresponding complex Lie algebra. The proof of this conjecture
was a major mathematical effort connecting these polynomials with
intersection cohomology of Schubert varieties, then using
Riemann-Hilbert correspondence to turn the intersection
cohomology complexes into $D$-modules on the flag variety, and
finally relating these to modules of the universal enveloping
algebra of the Lie algebra by global sections and localization
functors.

The last two steps in the proof of the conjecture were completed
by Beilinson-Bernstein \cite{BeilBern:KL} and Brylinski-Kashiwara
\cite{BrylKash:KL}. We will take this section to explain only the
first step, which was done by Kazhdan and Lusztig in
\cite{KaLus:Hecke2}. We follow the exposition of Springer
\cite{Springer:KLpols}, see also \cite{Soergel:ICpaper}.
Background on algebraic groups can be found for example in the
textbooks
\cite{Borel:AlgGroupBook,Humphreys:AlgGroupBook,Springer:AlgGroupBook}.

\subsection{Flag varieties and Schubert varieties}\label{s:FlagIntro}
Let $G$ be a reductive linear algebraic group over $\C$ with fixed
Borel subgroup $B$ and an algebraic torus $T\subset B$. The Weyl
group is $W=N_G(T)/T$. For $GL_n$ we may take $B$ to be the
subgroup of upper-triangular matrices and $T$ the diagonal
matrices. Then $W=S_n$. Let $\dot w\in N_G(T)$ be a
representative of $w\in W$, for example a permutation matrix in
$GL_n$. Denote by $S\subset W$ the set of simple reflections
generating $W$, so for example the adjacent transpositions in the symmetric group.
 The length $\ell(w)$ of $w\in W$ is the minimal
number of factors required to write $w$ as a product of simple
reflections.

The flag variety $X$ of $G$ is defined to be the homogeneous space
$G/B$. It is decomposed into  $B$-orbits by the {\it Bruhat
decomposition},
\begin{equation*}
G/B=\bigsqcup_{w\in W} B\dot w B/B,
\end{equation*}
where the individual orbits are affine cells $B\dot w
B/B\cong\C^{\ell(w)}$.

The Schubert varieties are defined to be the closures of the
Bruhat cells,
\begin{equation*}
X_w:=\overline{B\dot w B/B}.
\end{equation*}
The partial order on the Schubert varieties by inclusion gives
rise to a partial order on $W$ called the {\it Bruhat order}. It
is a general fact that all Borel subgroups are conjugate to one
another, and the normalizer of $B$ is just $B$. Therefore the
elements of the flag variety may be identified with Borel
subgroups via the correspondence
 $$
gB\leftrightarrow g\cdot B:=gBg\inv.
 $$

It will be useful to consider also the product $X\x X$. Let $ G$
act on $ X\x X$ diagonally. By an application of Bruhat
decomposition, the $G$-orbits in $X\x X$ are of the form
 $${\mathcal O}(w):=\text{ $
G$-orbit of $( B,\dot w\cdot B )$.}
 $$
We say that two elements $ B_1, B_2$  of the flag variety have
{\it relative position} $w$ if $( B_1, B_2)\in\mathcal O(w)$.
Write in this case $ B_1\overset w\to B_2$. The basic properties
of Bruhat decomposition imply
\begin{equation}\label{e:BruProps}
\begin{array}{lcl}
 B_1\overset{s}\to B_2\overset {s}\to  B_3 &\Rightarrow &
 B_1\overset{s}\to B_3\quad\text{or}\quad B_1=B_3,\\
 B_1\overset{v}\to  B_2\overset {w}\to  B_3 &\Rightarrow &
B_1\overset{vw}\to B_3 \quad
 \text{if}\quad \ell(vw)=\ell(v)+\ell(w),
\end{array}
\end{equation}
for $s\in S$ and $v,w\in W$. The set of $B'$ in $X$ such that
$B\overset w\to B'$ is just the Bruhat cell $B\dot w\cdot B$.

\subsection{The Hecke algebra} Let $v$ be
an indeterminate and let $\Hecke$ be a $\Z[v,v\inv]$-module with
basis indexed by $W$,
\begin{equation*}
\Hecke:=\Z[v,v\inv]\otimes\left<T_w\, |\, w\in W\right>_\Z.
\end{equation*}
There is a product on  $\Hecke$ with unit element $T_1$ such that
\begin{equation*}\label{e:Hrels}
\begin{array}{rlll}
T_s T_s&=(v^2-1)T_s+ v^2 T_1, &\quad &\text{for $s\in S$}\\
T_u T_w&= T_{uw}, & &\text{for $u,w\in W$ with $
\ell(uw)=\ell(u)+\ell(w)$.}
\end{array}
\end{equation*}
See for example the construction below, in Section~\ref{s:constr}. 
The algebra $\mathcal H$ determined by these relations 
is called the {\it Hecke algebra} associated to
the Weyl group $W$. When $v$ is specialized to $1$ it is clearly
just the group algebra $\Z[W]$ of $W$,

\subsection{Construction over $\mathbb F_q$}\label{s:constr}
There is a sense in which the Hecke algebra is the deformation of
$\Z[W]$ one obtains when replacing $W$ with the flag variety over
$\mathbb F_q$  (where $q$ corresponds to $v^2$). Let $\dot G,\dot
B, \dot X$ be split forms of $G,B$ and $X$ over the finite field
$\mathbb F_q$. For example $GL_n(\mathbb F_q)$. Then the
definitions from Section~\ref{s:FlagIntro} can be carried over to
this case.

Let us consider the convolution algebra $\Z[\dot X\x \dot
X]^{\dot G}$ of $\dot G$-invariant functions
on $\dot X\x \dot X$. 
Let $T_w\in \Z[\dot X\x \dot X]^{\dot G}$ be the characteristic
function of the $\dot G$-orbit $\dot{\mathcal O}(w)$. The
convolution product on $\Z[\dot X\x \dot X]^{\dot G}$ is
\begin{equation}\label{e:convolution}
S * T\, (\dot B_1,\dot B_3)=\sum_{\dot B_2\in \dot X}S(\dot
B_1,\dot B_2 ) T(\dot B_2,\dot B_3).
\end{equation}
So for example if $s\in S$
 $$
 T_s * T_s\, (\dot B_1,\dot B_3)=\begin{cases}
 q & \dot B_1=\dot B_3\\
 q-1 &\dot B_1\overset s\to \dot B_3\\
 0 & \text{otherwise.}
 \end{cases},
 $$
since there is a summand (of $1$) for every $\dot B_2$ with $\dot
B_1\overset s\to \dot B_2\overset s\to \dot B_3$. Therefore we
find
 $$
 T_s T_s=(q-1)T_s + q T_1.
  $$
The relation $T_u T_w= T_{uw}$ if $\ell(u)+\ell(w)=\ell(uw)$
follows in a similar way using \eqref{e:BruProps}.

\subsection{The Hecke algebra and perverse sheaves}
Above we have recovered the relations of the Hecke algebra from
the structure of the flag variety and its Bruhat decompositions
over $\mathbb F_q$. We can get a more sophisticated model by
replacing constructible  functions with semisimple perverse
sheaves. We return to working over $\C$. The stratification by
Bruhat cells of a Schubert variety $X_w$ and the $G$-orbit
stratification of $\overline{\mathcal O(w)}$ are both Whitney
stratifications, see
\cite[Proposition 1.4]{Dimca:SingularitiesBook}.

Let $\mathcal A^\bullet\in \mathcal D^b(X\x X)$ be constructible
for the $G$-orbit stratification of $X\x X$. Define
 $$
h_{X\x X}(\mathcal A^\bullet):=\sum_{w\in W}\sum_{i\in\Z}
\dim(\mathcal H^i_{(B,\dot w B)}(\mathcal A^\bullet))v^i T_w \in
\Hecke.
 $$
Similarly, for a bounded complex $\mathcal B^\bullet$ on $X$ that
is constructible for the Bruhat decomposition, define
 $$
h_{X}(\mathcal B^\bullet):=\sum_{w\in W}\sum_{i\in\Z}
\dim(\mathcal H^i_{\dot w B}(\mathcal B^\bullet))\, v^i T_w \in
\Hecke.
 $$
We have of course,
 $$
 h_{X\x X}(\Const_{\mathcal O(w),X\x X})=h_X(\Const_{B\dot w\cdot
 B,X})=T_w.
 $$
For another simple example let $w=s$ be a simple reflection. Then
we get that $X_s=(B\dot s\cdot B)\sqcup (1\cdot B )\cong \C P^1$
and therefore $\ICsheaf^\bullet_{X_s,X}=\Const_{X_s,X}[1]$. So
 $$
 h_X(\ICsheaf^\bullet_{X_s,X})=v^{-1}(T_s+T_1).
 $$

\begin{rem}
Consider the projection $p_1:\overline{\mathcal O(w)}\to X$ given
by $p_1(B_1,B_2)=B_1$. This is a bundle with fiber
$p_1\inv(B)=X_w$ (locally trivial by $G$-equivariance). Therefore
$\overline{\mathcal O(w)}$ looks locally like a product $\C^N\x
X_w$, where $N=\dim_\C(X)$. It follows that
$\ICsheaf^\bullet_{\overline{\mathcal O(w)}, X\x X}$ has stalks
given by
 $$
\left(\ICsheaf^\bullet_{\overline{\mathcal O(w)}, X\x X}\right
)_{(B,\dot v B)}= \left
(\ICsheaf^\bullet_{{X_w},X}[N]\right)_{\dot v B},
 $$
by a local application of Proposition~\ref{p:suspension}, for
example. Therefore
\begin{equation*}
h_{X\x X}(\ICsheaf^\bullet_{\overline{\mathcal O(w)},X\x X
})=v^{-N}h_X(\ICsheaf^\bullet_{X_w,X}).
\end{equation*}
\end{rem}

\subsection{Convolution product}
Suppose $\mathcal A^\bullet,\mathcal B^\bullet\in\mathcal D^b(X)$
are semisimple perverse sheaves constructible for the $G$-orbit
stratification. Consider the projections
\begin{equation*}
\begin{array}{ccccc}
X\x X &         &          &    & \\
      &\overset{p_{12}}\nwarrow&          &    & \\
      &         &X\x X\x X &\overset {p_{13}}\To & X\x X.\\
      &\overset{p_{23}}\swarrow &          &    & \\
X\x X &         &          &    & \\
\end{array}
\end{equation*}
Then
 $$
\mathcal A^\bullet\circ \mathcal B^\bullet
:=R(p_{13})_*(p_{12}^*\mathcal A^\bullet\otimes p_{23}^* \mathcal
B^\bullet)
 $$
is again a semisimple perverse sheaf by the decomposition theorem,
and constructible for the stratification by $G$-orbits.

For the purpose of simplifying notation let us write $w$ for the
element $\dot w\cdot B$ in $X$ or $(B,\dot w\cdot B)\in X\x X$,
when considering stalks at these points.

\begin{lem}\label{l:mult}
Suppose $\mathcal A^\bullet\in\mathcal D^b(X\x X)$ is
constructible with respect to $G$-orbits and has $\mathcal
H^i(\mathcal A^\bullet)=0$ for odd $i$. Then the same holds for
$\Const_{\overline{\mathcal O(s)},X\x X}\circ\mathcal A^\bullet$
and
 $$
 h_{X\x X}(\Const_{\overline{\mathcal O(s)},X\x X}\circ\mathcal
 A^\bullet)=(T_s+1)\cdot h_{X\x X}(\mathcal A^\bullet).
 $$
\end{lem}

\begin{proof}
The coefficient of $T_w$ in $(T_s+1)\cdot h_{X\x X}(\mathcal
A^\bullet)$ is
\begin{equation*}
\begin{cases}
\sum_{i\in\Z} v^{i}(\dim\mathcal H^i_w(\mathcal
A^\bullet)+\dim\mathcal H^{i-2}_{sw}(\mathcal A^\bullet))  & w<sw \\
\sum_{i\in\Z} v^{i}(\dim\mathcal H^i_{sw}(\mathcal A^\bullet)
+\dim \mathcal H^{i-2}_w(\mathcal A^\bullet)) & sw<w.
\end{cases}
\end{equation*}
We need to show, therefore, that
\begin{equation}\label{e:goal}
\dim (\mathcal H^i_w(\Const_{\overline{\mathcal O(s)},X\x X
}\circ\mathcal A^\bullet))=
\begin{cases} \dim\mathcal
H^i_w(\mathcal
A^\bullet)+\dim\mathcal H^{i-2}_{sw}(\mathcal A^\bullet)  & w<sw \\
\dim\mathcal H^i_{sw}(\mathcal A^\bullet) +\dim \mathcal
H^{i-2}_w(\mathcal A^\bullet) & sw<w.
\end{cases}
\end{equation}
Since $p_{13}$ is proper, $(p_{13})_*=(p_{13})_!$ and we have by
\eqref{e:ProperStalks} that
\begin{multline}\label{e:LHS}
\mathcal H^i_{w}(R(p_{13})_!(p_{12}^*\Const_{\overline {O(s)},X\x
X }\otimes p_{23}^*\mathcal A^\bullet))\\ = \mathbb
H^i_c(p_{13}\inv(B, \dot w\cdot B)\, ,\, p_{12}^*\Const_{\overline
{O(s)},X\x X}\otimes p_{23}^*\mathcal A^\bullet).
\end{multline}
Define
 $$D=\{(B_1, B_2, B_3 )\in p_{13}\inv(B,\dot w\cdot B)\ |\ (B_1,B_2)\in\overline{\mathcal
O(s)}\, \}.
 $$
Then $D\cong X_s \cong\C P^1$ and the right hand side in
\eqref{e:LHS} is simply $\mathbb H^i_c(D,p_{23}^*\mathcal
A^\bullet)$.

Let us now prove the identity \eqref{e:goal} for $sw<w$. In that
case $D$ decomposes into
\begin{align*}
D_1&=\{(B,B',\dot w\cdot B)\in D\ |\ B' \overset{w}\To \dot w\cdot B  \},\\
D_{0}&=\{(B,\dot s\cdot B,\dot w\cdot B)\},
\end{align*}
where $D_1\cong \C$. The inclusions
 $$
D_1\overset j\hookrightarrow D\overset{i}\hookleftarrow D_0
 $$
give rise to a distinguished triangle, $\,\cdots\to j_! j^!\to
1\to i_* i^*\overset{[1]}\to\cdots\, $, and a corresponding long
exact sequence
 $$
 \cdots\to\mathbb H^i_c(D_1,p_{23}^*\mathcal A^\bullet|_{D_1})\to
\mathbb H^i_c(D,p_{23}^*\mathcal A^\bullet)\to \mathbb
H^i_c(D_0,p_{23}^*\mathcal A^\bullet|_{D_0})\to\cdots.
 $$
Starting from the left, note that $p_{23}^*\mathcal A^\bullet$ has
constant cohomology  $\mathcal H^*_w(\mathcal A^\bullet)$ along
$D_1$. Therefore by Lemma~\ref{l:costalks} we have $\mathbb
H^i_c(D_1,p_{23}^*\mathcal A^\bullet|_{D_1})=\mathcal
H^{i-2}_w(\mathcal A^\bullet)$. Next
 $ \mathbb H^i_c(D,p_{23}^*\mathcal A^\bullet)=
 \mathcal H^i_w(\Const_{\overline{\mathcal O(s)},X\x X}\circ\mathcal
 A^\bullet)
 $
as we saw above. Finally, $\mathbb H^i_c(D_0,p_{23}^*\mathcal
A^\bullet|_{D_0})$ is trivially equal to $\mathcal
H^i_{sw}(\mathcal A^\bullet)$. So we have a long exact sequence
 $$
\cdots\to 0\to\mathcal H^{i-2}_w(\mathcal A^\bullet)\to \mathcal
H^i_w(\Const_{\overline{\mathcal O(s)},X\x X}\circ\mathcal
 A^\bullet) \to\mathcal H^i_{sw}(\mathcal
A^\bullet)\to 0\to\cdots ,
 $$
using that $\mathcal A^\bullet$ has vanishing cohomology in odd
degrees, and where $i$ is assumed even. This implies the equality
\eqref{e:goal} in the $sw<w$ case. The other case is similar.
\end{proof}
\begin{defn} Let $\iota:\Hecke\to\Hecke$ be the anti-involution
defined on generators by
\begin{equation*}
\begin{array}{lcl}
v &\mapsto& v\inv\\
T_s &\mapsto & T_{s}\inv=v^{-2} T_s + (v^{-2}-1).
\end{array}
\end{equation*}
Note that $\iota(T_s+1)=q^{-1}(T_s+1)$, where $q=v^2$.
\end{defn}

\begin{thm}[\cite{KaLus:Hecke}]\label{t:KL1} For any $w\in W$ there exists  a
unique element $C'_w\in\Hecke$ such that
\begin{enumerate}
\item $\iota (C'_w)=C'_w$.
\item The element $C'_w$ is expressed as
$$C'_w=v^{-{\ell(w)}}\sum_{u\le w}P_{u,w}(v^2) T_u$$ for
polynomials $P_{u,w}\in\Z[q]$ satisfying $\deg(P_{u,w})\le
\frac 12(\ell(w)-\ell(u)-1)$ if $u<w$ and such that $P_{w,w}=1$.
\end{enumerate}
\end{thm}
The polynomials $P_{u,w}$ are called the Kazhdan-Lusztig
polynomials. They can be computed by recursive formulas.

\begin{thm}\cite{KaLus:Hecke2}\label{t:KL2}
$C'_w=h_X(\ICsheaf^\bullet_{ X_w,X})$. In particular the
Kazhdan-Lusztig polynomials have nonnegative integer
coefficients.
\end{thm}

The remainder of these notes will be concerned with proving this
theorem. In other words, we want to show that
$h_X(\ICsheaf^\bullet_{X_w,X})$ satisfies the conditions in
Theorem~\ref{t:KL1}. The main step is to show that the involution
$\iota$ reflects the action of Verdier duality on the complexes of
sheaves.

Let $w$ be fixed, and with it a reduced expression
$w=s_{1}s_{2}\cdots s_{k}$ of $w$ in terms of the simple
reflections. Consider the variety
\begin{equation*}
Y=\{(B_1,\dotsc, B_{k+1})\in X^{k+1}\, |\,
(B_i,B_{i+1})\in\overline{\mathcal O(s_i)}\},
\end{equation*}
and its subvariety $Y_0=\{(B_1,\dotsc, B_{k+1})\in Y\, |\,
B_1=B\}$. Then we have projection maps
\begin{equation}\label{e:diag}
\begin{array} {ccccc}
Y_0&\hookrightarrow &Y &  \quad &        (B_1,\dotsc, B_{k+1})\\
\pi_0\downarrow & & \pi\downarrow & &           \downarrow   \\
X_w&\overset i\hookrightarrow &\overline{ \mathcal O(w)}& &
(B_1,B_{k+1}).
\end{array}
\end{equation}
Note that $Y$ and $Y_0$ are smooth. The map $\pi_0$ is the {\it
Bott-Samelson resolution} of the Schubert variety $X_w$.

\begin{lem}\label{l:assoc}
$R\pi_*(\Const_Y)=\Const_{\overline{\mathcal O(s_1)},X\x X}\circ
\Const_{\overline{\mathcal O(s_2)},X\x X}\circ\cdots \circ
\Const_{\overline{\mathcal O(s_k)},X\x X}.$
\end{lem}

Note that $\Const_Y=p_{12}^*\Const_{\overline{\mathcal O(s_1)},X\x
X}\otimes p_{23}^* \Const_{\overline{\mathcal O(s_2)},X\x
X}\otimes\cdots \otimes p_{k,k+1}^*\Const_{\overline{\mathcal
O(s_k)},X\x X}|_Y$ and apply $R\pi_*$.

\begin{proof}[Proof of Theorem~\ref{t:KL2}]
Firstly we need to prove that
$h_{X}(\ICsheaf^\bullet_{X_w,X})\in\Hecke$ is $\iota$-invariant.
If $w=1$ this is trivial. Also for  $s\in S$ the Schubert variety
$X_s\cong \C P^1$ and we know $\iota$-invariance by direct
calculation. For general $w$ the proof will be by induction. Let
us assume
$\iota(h_{X}(\ICsheaf^\bullet_{X_u,X}))=h_{X}(\ICsheaf^\bullet_{
X_u,X})$ for all $u<w$.

Consider the diagram \eqref{e:diag}. By proper base change and
the decomposition theorem,
\begin{equation*}
\closed^* R\pi_* (\Const_Y
[\ell(w)])=R{\pi_0}_*(\Const_{Y_0}[\ell(w)])= \bigoplus_{u\le w,
i\in\Z}\ICsheaf^\bullet_{X_u,X}\otimes \underline V_{u,i}[i]
\end{equation*}
for some multiplicity vector spaces $V_{u,i}$. The middle
expression is Verdier self-dual, since $\pi_0$ is proper and $Y_0$
smooth. On the other hand, the $\ICsheaf^\bullet_{X_u,X}$ on the
right hand side are each self-dual. Therefore we must have
 $$
 V_{u,i}\cong V_{u,-i}.
 $$
Also since $\pi_0|_{\pi\inv(X_w)}$ is an isomorphism,
\begin{equation*}
 V_{w,i}=\begin{cases}\C, & i=0,\\
 0, &i\ne 0.
 \end{cases}
\end{equation*}
Therefore
\begin{equation*}
h_{X}(\closed^* R\pi_* (\Const_Y[\ell(w)]))=
h_{X}(\ICsheaf^\bullet_{X_{w},X})+\sum_{u<w}h_{X}(\ICsheaf^\bullet_{X_u,X})p_u(v)
\end{equation*}
for Laurent polynomials $p_u$ with $p_u(v)=p_u(v\inv)$. Now by
Lemma~\ref{l:mult} and Lemma~\ref{l:assoc} we have
\begin{equation}\label{e:ing1}
h_X(\closed^* R\pi_*
(\Const_Y[\ell(w)]))=v^{-{\ell(w)}}(T_{s_1}+1)\cdots (T_{s_k}+1).
\end{equation}
Therefore it follows that $h_X(\closed^* R\pi_*
(\Const_Y[\ell(w)]))$ is invariant under the involution $\iota$.
Since
\begin{equation}\label{e:ing2}
h_{X}(\ICsheaf^\bullet_{X_w,X})=h_X(\closed^* R\pi_*
(\Const_Y[\ell(w)]))-\sum_{u<w}p_u(v)
h_{X}(\ICsheaf^\bullet_{X_u,X})
\end{equation}
it follows that $ h_{X}(\ICsheaf^\bullet_{X_w,X})$ is invariant
under $\iota\,$, by the induction hypothesis.

It remains to prove that $h_{X}(\ICsheaf^\bullet_{X_w,X})$
satisfies the second condition of Theorem~\ref{t:KL1}. Define
polynomials $\tilde P_{u,w}$ by
\begin{equation*}
h_{X}(\ICsheaf^\bullet_{X_w,X})=v^{-{\ell(w)}}\sum_{u\le w}\tilde
P_{u,w}(v) T_u.
\end{equation*}
Clearly, $\tilde P_{w,w}=1$. From the axiom $(C3)$ for the
Deligne sheaf on $X_w$ it follows that $\mathcal H^{-i}_u
(\ICsheaf^\bullet_{X_w,X})=0$ if $u<w$ and
$-i>-\ell(w)+2(\ell(w)-\ell(u))$.  This implies that $\deg(\tilde
P_{u,w})\le \ell(w)-\ell(u)-1$. Moreover,
\begin{equation*}
v^{{\ell(w)}} h_{X}(\ICsheaf^\bullet_{X_w,X})\in\Z[v^2],
\end{equation*}
by equations \eqref{e:ing1} and \eqref{e:ing2} and induction. So
$\tilde P_{u,w}(v)=P_{u,w}(v^2)$ where the $P_{u,w}$ are
polynomials satisfying the axioms for the Kazhdan-Lusztig
polynomials.
\end{proof}
\bibliographystyle{amsplain}


\def\cprime{$'$}
\providecommand{\bysame}{\leavevmode\hbox to3em{\hrulefill}\thinspace}
\providecommand{\MR}{\relax\ifhmode\unskip\space\fi MR }
\providecommand{\MRhref}[2]{%
  \href{http://www.ams.org/mathscinet-getitem?mr=#1}{#2}
}
\providecommand{\href}[2]{#2}

\end{document}